\documentclass{amsart}
\usepackage[latin1]{inputenc}
\usepackage{amssymb}
\usepackage{amsmath}
\usepackage{amsbsy}
\usepackage{enumerate}
\usepackage{epsfig,color,caption}
\usepackage[dvipsnames]{xcolor}
\usepackage{colonequals}
\usepackage{makecell}
\usepackage{color}
\usepackage[pdftex,colorlinks,citecolor=blue]{hyperref}
\usepackage{cleveref}
\usepackage{tikz}
\usetikzlibrary{arrows}
\usetikzlibrary{positioning}
\usetikzlibrary{cd}
\usetikzlibrary{calc}
\usetikzlibrary{shapes}
\usepackage{comment}
\usepackage{longtable,graphics,multirow,ulem}
\usepackage{changepage}
\usetikzlibrary{matrix}
\usepackage{pgf,tikz}
\usepackage[all]{xy}
\newtheorem{theorem}{Theorem}[section]
\newtheorem{lemma}[theorem]{Lemma}
\newtheorem{proposition}[theorem]{Proposition}
\newtheorem{corollary}[theorem]{Corollary}

\newtheorem{definition}[theorem]{Definition}

\newtheorem{rem}[theorem]{Remark}

\numberwithin{equation}{section}
\newcommand{\rk}{\mbox{rank }}
\newcommand{\NS}{\mbox{NS}}

\newcommand{\ra}{\rightarrow}

\newcommand{\C }{ \mathbb{C}}

\newcommand{\Z}{\mathbb{Z}}
\newcommand{\Q}{\mathbb{Q}}

\makeatletter
\def\blfootnote{\xdef\@thefnmark{}\@footnotetext}

\title[Order 3 symplectic automorphisms on K3 surfaces]{Order 3 symplectic automorphisms on K3 surfaces}
\begin{document}

\author{Alice Garbagnati}
\address{Alice Garbagnati, Dipartimento di Matematica, Universit\`a Statale di Milano,
	via Saldini 50, I-20133 Milano, Italy}
\email{alice.garbagnati@unimi.it}

\author{Yulieth Prieto Monta\~{n}ez}
\address{Yulieth Prieto Monta\~{n}ez, Dipartimento di Matematica,  Universit\`a di Bologna, Piazza di Porta S. Donato, 5, I-40126 Bologna, Italy}
\email{yulieth.prieto2@unibo.it}

\thanks{This paper was written during 2020, it cites and and greatly benefits of works by Conway, who died of Covid. We would like to remember him and all the victims of the pandemic. At this time (February 2021) it is estimated that they are 2209195}

\begin{abstract}
	
	The aim of this paper is to generalize results known for the symplectic involutions on K3 surfaces to the order 3 symplectic automorphisms on K3 surfaces. In particular, we will explicitly describe the action induced on the lattice $\Lambda_{K3}$, isometric to the second cohomology group of a K3 surface, by a symplectic automorphism of order 3; we exhibit the maps $\pi_*$ and $\pi^*$ induced in cohomology by the rational quotient map $\pi:X\dashrightarrow Y$, where $X$ is a K3 surface admitting an order 3 symplectic automorphism $\sigma$ and $Y$ is the minimal resolution of the quotient $X/\langle\sigma\rangle$; we deduce the relation between the N\'eron--Severi group of $X$ and the one of $Y$. Applying these results we describe explicit geometric examples and generalize the Shioda--Inose structures.
\end{abstract}	
\maketitle

\section{Introduction}
In his seminal paper \cite{Nik}, Nikulin introduced the finite order symplectic automorphisms of K3 surfaces and studied their properties considering  both their action on the surface (in particular determining their fixed locus) and the action they induce on the second cohomology group of the K3 surfaces. He proved that both the topology of the fixed locus and the action induced in cohomology are ``unique" and depend only on the order of the automorphism (see \cite[Definition 4.6, Theorem 4.7, Section 5]{Nik} for the precise definitions and results). 

One of the main reasons of interest in these particular automorphisms is that their existence creates a relation between two different (families of) K3 surfaces: the family of K3 surfaces admitting a symplectic automorphism of order $n$ and the one of the desingularization of the quotient of the K3 surfaces by a symplectic automorphism of order $n$. The latter surfaces are still K3 surfaces, but in general not isomorphic to the original ones. 

After \cite{Nik}, \cite{GS}, \cite{GS2} it is possible to describe the family $\mathcal{S}$ (resp. $\mathcal{T}$) of the projective K3 surfaces $X$ admitting a symplectic automorphism $\sigma_n$ of order $n$ for prime numbers $n$ (resp. of the projective K3 surfaces $Y$ obtained as desingularization of the quotient of a K3 surface by a symplectic automorphism of order $n$), in terms of families of the lattice polarized K3 surfaces. Both these families have countable many irreducible components. It remains an open problem to determine the relationship between the components of these two families. More explicitly, if one considers a K3 surface $X$ admitting a symplectic automorphism $\sigma_n$ of order $n$, $X$ is contained in a specific component of $\mathcal{S}$. One wants to use this information to determine the component of $\mathcal{T}$ which contains $Y$ (resolution of $X/\langle\sigma_n\rangle$). By the theory of the moduli space of lattice polarized K3 surfaces, see e.g. \cite{Do}, this problem is equivalent to determine the relation between the N\'eron--Severi group of $X$ and the one of $Y$. 

This problem has an answer for symplectic involutions, thanks to the results contained in \cite{Mo}, \cite{vGS} and \cite{GS2}.

More in general, the order 2 case is the most studied and several nice results are based on the following result by D. Morrison. In \cite{Mo}, the isometry induced by a symplectic involution on the lattice $\Lambda_{K3}:=U^{\oplus 3}\oplus E_8^{\oplus 2}$, the second cohomology group of any K3 surface, is described: it acts as the identity on the direct sum $U^{\oplus 3}$ and switches the two copies of $E_8$ in the summand $E_8^{\oplus 2}$.  
Thanks to this result Morrison finds a nice and strong relation, called Shioda--Inose structure, between certain K3 surfaces and an associated Abelian surface.

By using Morrison's result, van Geemen and Sarti describe the family $\mathcal{S}$ and the maps $\pi^*$ and $\pi_*$, where $\pi$ is the rational map $\pi:X\dashrightarrow Y$ induced by the quotient map $X\ra X/\langle \sigma_2\rangle$. As a consequence, in \cite{GS2} the description of the family $\mathcal{T}$ and the explicit relation between components of $\mathcal{S}$ and of $\mathcal{T}$ are obtained. This result is applied in \cite{CG} to construct infinite towers of isogenous K3 surfaces.

The aim of this paper is to generalize these results, known for an automorphism of order 2, to the order 3 case. As mentioned above, the results in the order 2 case mainly depend on the description of the isometry induced by a symplectic involution on the second cohomology group. Hence our first goal is to establish a similar result for the order 3 case. It turns out that the action induced by a symplectic automorphism of order 3 on $\Lambda_{K3}$ is not compatible with the direct summands $U^{\oplus 3}\oplus E_8^{\oplus 2}$. So we give a different description of the lattice $\Lambda_{K3}$ and we provide the action on such a lattice (see Theorem \ref{prop: NS(S), H2(S), action of sigma} and Section \ref{subsec: isometry on H2(S,Z)}):

{\bf Theorem A} {\it The lattice $\Lambda_{K3}$ is an overlattice of index $3^2$ of the lattice $A_2(-1)\oplus U\oplus E_6^{\oplus 3}$. The action induced by an order 3 symplectic automorphism on  $A_2(-1)\oplus U\oplus E_6^{\oplus 3}$ is the identity on the summands $A_2(-1)\oplus U$ and a cyclic permutation of the three summands in $E_6^{\oplus 3}$. The action on the overlattice $\Lambda_{K3}$ is induced by the one on $A_2(-1)\oplus U\oplus E_6^{\oplus 3}$ by $\Q$-linear extension.}

This result is achieved in a geometric way: thanks to Nikulin's theorem the isometry induced on $\Lambda_{K3}$ by a specific symplectic automorphism $\sigma$ of order 3 on a specific K3 surface $X$ is essentially unique (i.e. it does not depend on $X$ and $\sigma$). So in Section \ref{sec: a particular surface} we describe this action for a particular K3 surface $S$ with a specific automorphism $\sigma$ of order 3 and this suffices to prove the general result. The surface $S$ is a classic K3 surface (studied e.g. in \cite{SI} and \cite{V}), with Picard number 20. The order 3 automorphism $\sigma$ acts on the classes generating the N\'eron--Severi group of $S$ permuting three copies of $E_6$, see Theorem \ref{prop: NS(S), H2(S), action of sigma}. This allows us to describe the isometry induced by any automorphism of order 3 on the lattice $\Lambda_{K3}$, see Section \ref{subsec: isometry on H2(S,Z)}.

In Section \ref{sec: cohomological action} we describe the maps $\pi_*$ and $\pi^*$ induced by the quotient map.

In Section \ref{sec: projective families} we describe the families $\mathcal{S}$ and $\mathcal{T}$, listing all their components (Theorems \ref{theor:NS X} and \ref{theor:NS Y} respectively).

In Section \ref{sec: relations between families}, merging the results on the families $\mathcal{S}$  and $\mathcal{T}$ with the ones about the maps $\pi^*$ and $\pi_*$, we deduce the relations between the components of $\mathcal{S}$ and the ones of $\mathcal{T}$ and in particular we relate the N\'eron--Severi group of a projective K3 surface $X$ admitting a symplectic automorphism $\sigma$ of order 3 with the one of the desingularization $Y$ of $X/\langle\sigma\rangle$. This result is proved under assumptions of generality for $X$ and $Y$, it is stated in Theorem \ref{theo: NS(Y) iff NS(X)} and it can be resumed in the following: 

{\bf Theorem B} {\it Let $X$ be a projective K3 surface admitting an order 3 symplectic automorphism $\sigma$, then $\rho(X)\geq 13$ and if $\rho(X)=13$ the N\'eron--Severi group of $X$ determines uniquely the one of the minimal resolution of $X/<\sigma>$ and viceversa.}\\

The application of these results to specific families of K3 surfaces allows us to exhibit equations of families of projective K3 surfaces admitting symplectic automorphisms of order 3 and of their quotients in Section \ref{sec: examples}. For example we show that there is exactly one component of $\mathcal{S}$ which corresponds to family of quartic hypersurfaces in $\mathbb{P}^3$ admitting a symplectic automorphism of order 3 (and in this case the quotient is a singular surface in $\mathbb{P}^{7}$). But there are two different components of $\mathcal{S}$ which correspond to two families of complete intersections of type $(2,3)$ in $\mathbb{P}^4$ admitting a symplectic automorphism of order 3. The quotient of the surfaces contained in one of these components is a singular double cover of $\mathbb{P}^2$, and it corresponds to a component in $\mathcal{T}$, the quotient of the surfaces in the other case has a singular model in $\mathbb{P}^{10}$ and it corresponds to another component in $\mathcal{T}$, see Section \ref{subsec: d=3}.

In Section \ref{sec: applications} we briefly present two applications of the theorems A and B: thanks to theorem B one is able to construct infinite towers of isogenous K3 surfaces by considering iterated quotients by symplectic automorphisms of order 3. The analogous construction which uses symplectic involutions was presented in \cite{CG}. Thanks to Theorem A one is able to generalize the definition and previous results on Shioda--Inose structures, which are particular relations between an Abelian surface and a K3 surface. The analogous construction involving involutions was presented in \cite{Mo}.


{\bf Notation:} In all the paper the $ADE$ lattices are assumed to be negative definite.

Given a lattice $L$ we denote by $L(n)$ the lattice obtained from $L$ multiplying the bilinear form by $n$.

The lattice $L^*(n)$ is the scaled dual lattice of $L$, i.e. it is the dual of $L$ with the induced form multiplied by $n$. In particular $A_2^*(3)$ is an even negative definite lattice of rank 2 and discriminant group $\Z/3\Z$ and $E_6^*(3)$  is an even negative definite lattice of rank 6 and discriminant group $(\Z/3\Z)^5$.

The lattice $U$ is the unique even indefinite unimodular lattice of rank 2.

The lattice $K_{12}$ is the opposite of the Coxeter--Todd lattice, defined in \cite{CT}. It is an even negative definite lattice of rank 12 and discriminant group $(\Z/3\Z)^6$.\\

{\it Acknowledgement: We warmly thank Bert van Geemen and Giovanni Mongardi who read a preliminary version of this paper. We are grateful to the referee for the valuable suggestions.
The second author is supported by the fellowship INDAM-DP-COFUND-2015 "INdAM Doctoral Programme in Mathematics and/or Applications Cofunded by Marie Sklodowska-Curie Actions", Grant Number 713485.}


\section{An example of a K3 surface with a symplectic automorphism of order 3}\label{sec: a particular surface}
The aim of this section is to describe a specific K3 surface $S$ and a symplectic automorphism  $\sigma$ of order 3 on $S$. We will describe the geometry of $S$ and of its automorphism and we determine the action of $\sigma^*$ on the N\'eron--Severi group $\NS(S)$. From this description one deduces the action of $\sigma^*$ on $H^2(S,\mathbb{Z})\simeq \Lambda_{K3}$. 

\subsection{The surface $S$ and the automorphism $\sigma$}\label{subsec: the surface S and the automorphism sigma}

	Let $\zeta_3$ be a third root of unity and $E_{\zeta_3}$ the elliptic curve admitting an endomorphism of order 3, denoted by $\alpha$.
Let $S$ be the unique K3 surface obtained as minimal resolution of the quotient  $\left(E_{\zeta_3}\times E_{\zeta_3}\right)/\langle\alpha\times \alpha^2\rangle$. This surface is well known and it is one of the so called two most algebraic K3 surfaces, see e.g. \cite{V} and \cite{SI}. In particular it is known that the transcendental lattice of $S$ is $T_S\simeq A_2(-1)$ and hence the rank of the N\'eron--Severi group of $S$ is 20.

The surface $S$ is naturally endowed with an elliptic  fibration, induced by the trivial elliptic fibrations on $E_{\zeta_3}\times E_{\zeta_3}$. To find a Weierstrass equation of this fibration, we denote by $v_i^2=u_i^3-1$ the equation of the $i$-th copy of $E_{\zeta_3}$ in $E_{\zeta_3}\times E_{\zeta_3}$ and we assume the automorphism $\alpha$ to be $(v_i,u_i)\mapsto(v_i,\zeta_3 u_i)$. Hence the functions $y=v_1u_2^6$, $x=u_1u_2^4$, $v_2$ are invariant for the action of $\alpha\times \alpha^2$ and satisfy the following equation:
$$y^2=x^3-(v_2^2+1)^4$$ whose homogeneous form is 
$$y^2=x^3-(v_2^2+w_2^2)^4w_2^4.$$
Since the map induced by $y=v_1u_2^6$, $x=u_1u_2^4$, is a $3:1$ map $(E_{\zeta_3}\times E_{\zeta_3})\ra \left(E_{\zeta_3}\times E_{\zeta_3}\right)/\langle\alpha\times \alpha^2\rangle$, the previous equation is a Weierstrass equation of the surface $\left(E_{\zeta_3}\times E_{\zeta_3}\right)/\langle\alpha\times \alpha^2\rangle$ (cf. e.g. \cite[Section 7.1]{GP} to more details on the computation of equations of quotient surfaces).
By applying the change of coordinates $t:=v_2+\sqrt{3}w_2$, $s:=v_2-\sqrt{3}w_2$ (and multiplying the last term by $6\sqrt{3}$) one obtains the following Weierstrass equation for $S$:
\begin{eqnarray}\label{eq X}y^2=x^3+(t^3-s^3)^4.\end{eqnarray}

Let us now define the order 3 symplectic automorphism $\sigma$ of $S$. By \eqref{eq X}, it is clear that the fibration admits the sections $t\mapsto (x(t), y(t))=(0, (t^3-s^3)^2)$ and $t\mapsto (x(t), y(t))=(0, -(t^3-s^3)^2)$. One can directly check that these are 3-torsion sections, denoted in the following by $T_1$ and $T_2$. So, denoted by $\mathcal{O}$ the zero section of the fibration, one obtains that $\{\mathcal{O}, T_1, T_2\}$ with the group structure of the Mordell--Weil group is a copy of $\mathbb{Z}/3\mathbb{Z}$.
In particular, the translation by $T_1$ is a well defined automorphism of $S$ of order 3. This is the automorphism $\sigma$. 
It acts as a translation on each fiber of the fibration and as the identity on the basis, hence it is a symplectic automorphism on $S$.
\subsection{The N\'eron--Severi group of $S$ and the isometry $\sigma^*$ on $\NS(S)$.}\label{subsec: NS(S) and action of sigma}

The elliptic fibration \eqref{eq X} has three reducible fibers of type $IV^*$, over $t=1,\pm\zeta_3$, see e.g. \cite{SI} or \cite[Table IV.3.1]{Mi}. We recall that the fibers of type $IV^*$ are the ones whose components intersect as in the diagram $\widetilde{E_6}$, the extended Dynkin-diagram of $E_6$.
We denote by $C_i^{(j)}$ the $i$-th component of the $j$-th reducible fiber, with the following conventions:
$$\begin{array}{c}C_0^{(j)}C_1^{(j)}=C_1^{(j)}C_2^{(j)}=C_2^{(j)}C_3^{(j)}=C_3^{(j)}C_4^{(j)}=C_5^{(j)}C_6^{(j)}=C_5^{(j)}C_2^{(j)}=1\\ \\
	C_0^{(j)}\mathcal{O}=C_4^{(j)}T_1=C_6^{(j)}T_2=1, \ \ \left(C_i^{(j)}\right)^2=\mathcal{O}^2=T_1^2=T_2^2=-2\end{array}$$
and all the other intersections are 0. 
In terms of the classes $F$ (the class of the fiber), $\mathcal{O}$ and $C_i^{(j)}$, $i=1,\ldots, 6$, $j=1,2,3$, we have 
$$\begin{array}{c}T_1=2F+\mathcal{O}-\frac{1}{3}\sum_{j=1}^3\left(3C_1^{(j)}+6C_2^{(j)}+5C_3^{(j)}+4C_4^{(j)}+4C_5^{(j)}+2C_6^{(j)}\right),\\
	\\
T_2=2F+\mathcal{O}-\frac{1}{3}\sum_{j=1}^3\left(3C_1^{(j)}+6C_2^{(j)}+4C_3^{(j)}+2C_4^{(j)}+5C_5^{(j)}+4C_6^{(j)}\right),\\
\\

F=C_0^{(j)}+2C_1^{(j)}+3C_2^{(j)}+2C_3^{(j)}+C_4^{(j)}+2C_5^{(j)}+C_6^{(j)},\ \ j=1,2,3.\end{array}$$

Since $\sigma$ preserves each fiber and acts on each fiber as a translation by $T_1$, its action on the classes described above is the following:
$$\begin{array}{ccc}\mathcal{O}\stackrel{\sigma^*}{\mapsto} T_1\stackrel{\sigma^*}{\mapsto} T_2\stackrel{\sigma^*}{\mapsto}\mathcal{O}&&C_0^{(j)}\stackrel{\sigma^*}{\mapsto} C_4^{(j)}\stackrel{\sigma^*}{\mapsto} C_6^{(j)}\stackrel{\sigma^*}{\mapsto}C_0^{(j)}\\
C_1^{(j)}\stackrel{\sigma^*}{\mapsto} C_3^{(j)}\stackrel{\sigma^*}{\mapsto} C_5^{(j)}\stackrel{\sigma^*}{\mapsto}C_1^{(j)}&&
C_2^{(j)}\stackrel{\sigma^*}{\mapsto} C_2^{(j)}\stackrel{\sigma^*}{\mapsto} C_2^{(j)}\stackrel{\sigma^*}{\mapsto}C_2^{(j)}
\end{array}$$
where we used that for each pair of divisors $D_1$ and $D_2$ in $\NS(S)$, $D_1D_2=\sigma^*(D_1)\sigma^*(D_2)$.
The following three orthogonal copies of $E_6$ are permuted by $\sigma^*$:
\begin{align}\label{eq: 3E6}\xymatrix{C_1^{(1)}\ar@{-}[r]&C_0^{(1)}\ar@{-}[r]&\mathcal{O}\ar@{-}[r]\ar@{-}[d]&C_0^{(2)}\ar@{-}[r]&C_1^{(2)}\\&&C_0^{(3)}\\
	C_3^{(1)}\ar@{-}[r]&C_4^{(1)}\ar@{-}[r]&T_1\ar@{-}[r]\ar@{-}[d]&C_4^{(2)}\ar@{-}[r]&C_3^{(2)}\\&&C_4^{(3)}\\C_5^{(1)}\ar@{-}[r]&C_6^{(1)}\ar@{-}[r]&T_2\ar@{-}[r]\ar@{-}[d]&C_6^{(2)}\ar@{-}[r]&C_5^{(2)}\\&&C_6^{(3)}}\end{align}

The components of these three copies of $E_6$ are 18 classes in $\NS(S)$, which are independent, since their intersection form is $E_6^{\oplus 3}$.

Since $\rk(\NS(S))=20$, there are other two independent classes: the classes  $C_2^{(3)}$ and $D:=3\mathcal{O}+\left(C_1^{(1)}+C_1^{(2)}+C_1^{(3)}+2\left(C_0^{(1)}+C_0^{(2)}+C_0^{(3)}\right)\right)$  are orthogonal to the 3 copies of $E_6$ and their intersection form is $\left[\begin{array}{cc}-2&1\\1&0\end{array}\right]\simeq U$.
So the classes appearing in \eqref{eq: 3E6}, $C_2^{(3)}$ and $D$ are 20 linearly independent classes on $\NS(S)$ and thus they form a $\mathbb{Q}$-basis of $\NS(S)$.

The lattice generated by these classes is $U\oplus E_6\oplus E_6\oplus E_6$. Moreover there is an extra class in $\NS(S)$ given by $$G:=\frac{C_1^{(1)}+2C_0^{(1)}+C_0^{(2)}+2C_1^{(2)}+C_3^{(1)}+2C_4^{(1)}+2C_3^{(2)}+C_4^{(2)}+C_5^{(1)}+2C_6^{(1)}+2C_5^{(2)}+C_6^{(2)}}{3}.$$
Indeed, recalling that $C_0^{(j)}=F-2C_1^{(j)}-3C_2^{(j)}-2C_3^{(j)}-C_4^{(j)}-2C_5^{(j)}-C_6^{(j)},$
one obtains that  
\begin{eqnarray}\label{eq: G, the class x in NS(S)}%
G=F-C_1^{(1)}-2C_2^{(1)}-C_3^{(1)}-C_5^{(1)}-C_2^{(2)}.
\end{eqnarray}
Therefore $G$ is a linear combination with integer coefficients of the classes contained in $\NS(S)$, and thus it is contained in $\NS(S)$. 

So there is an overlattice of index 3 of $E_6^{\oplus 3}$ contained in $\NS(S)$, which is generated by the three copies of $E_6$ in \eqref{eq: 3E6} and by the class $G$. Let us denote this lattice by $(E_6^{\oplus 3})'$.
Then $\NS(S)\simeq U\oplus (E_6^{\oplus 3})'$, indeed we proved that $U\oplus (E_6^{\oplus 3})'\subset  \NS(S)$ since we exhibits the classes of $\NS(S)$ which generate $U\oplus (E_6^{\oplus 3})'$. The inclusion $U\oplus (E_6^{\oplus 3})'\subset  \NS(S)$ has a finite index, indeed the ranks of the two lattices are the same, and the index of the inclusion is 1, because the discriminants of these two lattices are the same.

The action of $\sigma^*$ on $U\oplus E_6^{\oplus 3}$ is very easy to describe: it acts as the identity on $U$ and permutes the three copies of $E_6$. Moreover $\sigma^*$ preserves the divisor $G$, since it is symmetric in the three copies of $E_6$.

\subsection{The lattice $H^2(S,\Z)$ and the isometry $\sigma^*$ on $H^2(S,\Z)$}\label{subsec: isometry on H2(S,Z)}

Since each symplectic automorphism acts as the identity on the generator of $H^{2,0}(S)$, it acts as the identity on the transcendental lattice $T_S$, i.e.
 $\sigma^*_{|T_S}=id_{T_S}$. We recall that $T_S\simeq A_2(-1)$.
So $H^2(S,\mathbb{Z})$ is an overlattice (of index 3) of $T_S\oplus \NS(S)$ and $\NS(S)$ is an overlattice of index 3 of $U\oplus E_6^{\oplus 3}$.
Hence we obtain that $H^2(S,\mathbb{Z})$ is an overlattice of index $3^2$ of $A_2(-1)\oplus U\oplus E_6\oplus E_6\oplus E_6$.
The action induced by $\sigma^*$ on $A_2(-1)\oplus U\oplus E_6^{\oplus 3}$ is trivial on $A_2(-1)\simeq T_S$, is trivial on the sublattice $U\subset \NS(S)$ generated by $C_2^{(3)}$ and $D$, and permutes the three copies of $E_6$ described in \eqref{eq: 3E6}, so the action of $\sigma^*$ is the following
\begin{equation}\label{eq:sigma* on H2}\begin{array}{ccccccccccccccccccccccc}\sigma^*:&A_2(-1)\oplus &U\oplus&E_6\oplus&E_6\oplus&E_6&\ra&A_2(-1)\oplus &U\oplus&E_6\oplus&E_6\oplus&E_6\\
&(a,&u,&e,&f,&g)&\mapsto&(a,&u,&g,&e,&f)\end{array}\end{equation}

Since $H^2(S,\mathbb{Z})$ is an overlattice of finite index of $A_2(-1)\oplus U\oplus E_6^{\oplus 3}$, the action of $\sigma^*$ on $H^2(S,\mathbb{Z})$ is induced by the $\mathbb{Q}$-linear extension of \eqref{eq:sigma* on H2}. 

In order to obtain $H^2(S,\mathbb{Z})$ from $A_2(-1)\oplus U\oplus E_6^{\oplus 3}$ one has to add two classes, which are contained in the discriminant group of $A_2(-1)\oplus U\oplus E_6^{\oplus 3}$.

Let us fix the following generators for the discriminant group of $E_6$ and $A_2(-1)$ respectively:
\begin{equation}\label{eq: element in discriminant group}v^{(i)}=\frac{e_1^{(i)}+2e_2^{(i)}+e_4^{(i)}+2e_5^{(i)}}{3},\ i=1,2,3\mbox{ and }w=\frac{a_1+2a_2}{3}\end{equation}
where $a_h$ are the generators of $A_2(-1)$ and $e_j^{(i)}$, $j=1,\ldots 6$ $i=1,2,3$ are the generators of $E_6$ whose intersections are as in the diagram:
\begin{eqnarray}\label{eq:E6}\xymatrix{e_1^{(i)}\ar@{-}[r]&e_2^{(i)}\ar@{-}[r]&e_3^{(i)}\ar@{-}[r]\ar@{-}[d]&e_4^{(i)}\ar@{-}[r]&e_5^{(i)}\\&&e_6^{(i)}}\end{eqnarray}

We consider the overlattice of $E_6^{\oplus 3}$ obtained by adding to $E_6^{\oplus 3}$ the class $$x:=v^{(1)}+v^{(2)}+v^{(3)}.$$ It coincides with the lattice $(E_6^{\oplus 3})'$ constructed before (an explicit isometry is obtained mapping the generators of $E_6^{\oplus 3}$ to the curves which appear in \eqref{eq: 3E6} and the vector $x$ to the class $G$ defined in \eqref{eq: G, the class x in NS(S)}).
 
Hence we constructed the lattice $A_2(-1)\oplus U\oplus (E_6^{\oplus 3})'$ as an overlattice of index 3 of  $A_2(-1)\oplus U\oplus E_6^{\oplus 3}$ by adding the class $x$.

Since $\left(v^{(i)}\right)^2=-4/3$ and $w^2=2/3$, the class $y=w+v^{(1)}-v^{(2)}$ is such that $y^2=-2$ (equivalent to 0 mod $2\Z$). In particular $y\in A_{A_2(-1)\oplus U\oplus (E_6^{\oplus 3})'}$ is isotropic. Hence by adding $y$ to $A_2(-1)\oplus U\oplus (E_6^{\oplus 3})'$ one obtains an overlattice of $A_2(-1)\oplus U\oplus (E_6^{\oplus 3})'$ which is even, unimodular and whose signature is $(3,19)$. This implies then this lattice is $\Lambda_{K3}$. Since $A_2(-1)\simeq T_S$ and $U\oplus (E_6^{\oplus 3})'\simeq \NS(S)$, $y$ is the gluing vector, which enlarges $T_S\oplus \NS(S)$ to $H^2(S,\mathbb{Z})$. 


To resume the results of this section we proved the following (with the previous notations): 
\begin{theorem}\label{prop: NS(S), H2(S), action of sigma}
	There is an even overlattice $(E_6^{\oplus 3})'$ of index 3 of $E_6^{\oplus 3}$ which is obtained by adding to $E_6^{\oplus 3}$ the class $x=v^{(1)}+v^{(2)}+v^{(3)}$ where $v^{(i)}$ are as in \eqref{eq: element in discriminant group}.
	
	The N\'eron--Severi group of the K3 surface $S$ is isometric to $U\oplus \left(E_6^{\oplus 3}\right)'$ and the action of $\sigma^*$ on  $U\oplus E_6^{\oplus 3}$ is a cyclic permutation of order 3 on $E_6^{\oplus 3}$ and is the identity on $U$.
	
	The lattice $H^2(S,\Z)$ is an overlattice of index $3^2$ of $A_2(-1)\oplus U\oplus E_6^{\oplus 3}$ obtained by adding $x$ and $y:=w+v^{(1)}-v^{(2)}$ where $w$ and $v^{(i)}$ are as in \eqref{eq: element in discriminant group}. 
	
	The action of $\sigma^*$ on $A_2(-1)\oplus U\oplus E_6^{\oplus 3}$ is the identity on $A_2(-1)\oplus U$ and a cyclic permutation of order 3 on $E_6^{\oplus 3}$; the one on $H^2(S,\Z)$ is obtained extending this one to $x$ and $y$.
\end{theorem}



\section{The cohomological action induced by a symplectic automorphism of order 3}\label{sec: cohomological action}
In the previous Section we described the action of a particular symplectic automorphism of order 3 on a particular K3 surface and the isometry it induces.
Since the action of a finite order symplectic automorphism on the second cohomology group of a K3 surface is unique (cf. \cite[Definition 4.6 and Theorem 4.7]{Nik}), one obtains general results by the ones described above in a specific example. The aim of this section is to state and prove these more general results: first, in Section \ref{subsect: the action on Lambda} we describe the isometry induced by a symplectic automorphism of order 3 on the lattice $\Lambda_{K3}$; then, in Sections \ref{subsec: the map pi_*} and \ref{subsec: pi*} we describe the maps $\pi_*$ and $\pi^*$, where $\pi$ is induced by the quotient map. The description of these maps is the main technical result of our paper and it is the analoguos of the results proved in \cite{vGS} for the order 2 case.

Some properties of the action of $\sigma^*$ on $\Lambda_{K3}$ were already known. In particular the lattice $\left(H^2(S,\Z)^{\sigma^*}\right)^{\perp}$ is isometric to the lattice $K_{12}$, by \cite[Theorem 4.1]{GS}. The lattice $K_{12}$ was defined by Coxeter and Todd, in \cite{CT} and intensively studied in \cite{CS} by Conway and Sloane where its relation with the lattice $E_6$ is analysed from a different point of view.

\subsection{The action induced by an automorphism of order 3 on $\Lambda_{K3}$}\label{subsect: the action on Lambda}
Recall that $\Lambda_{K3}\simeq U^{\oplus 3}\oplus E_8^{\oplus 2}$ is the unique even unimodular lattice of signature $(3,19)$.
Since for each K3 surface $X$, $H^2(X,\mathbb{Z})$ is an even unimodular lattice of signature $(3,19)$, $H^2(X,\mathbb{Z})$ is isometric to $\Lambda_{K3}$. By Theorem \ref{prop: NS(S), H2(S), action of sigma} we have an alternative description of $\Lambda_{K3}$: it is the overlattice of index $3^2$ of $A_2(-1)\oplus U\oplus E_6^{\oplus 3}$ obtained by adding the classes $x$ and $y$.

Let $\sigma$ be a symplectic automorphism of order 3 on a K3 surface $X$, then $\sigma^*$ acts on $H^2(X,\mathbb{Z})$ and this action is unique, by \cite[Theorem 4.7]{Nik}. So $\sigma^*$ is an order 3 isometry of $\Lambda_{K3}$. In Theorem \ref{prop: NS(S), H2(S), action of sigma} we described this isometry by considering a very special K3 surface and by the uniqueness of this action we conclude that
\begin{align}\label{eq: action sigma}\begin{array}{ccccccccccccccccccccccc}\sigma^*:&A_2(-1)\oplus &U\oplus&E_6\oplus&E_6\oplus&E_6&\ra&A_2(-1)\oplus &U\oplus&E_6\oplus&E_6\oplus&E_6\\
&(a,&u,&e,&f,&g)&\mapsto&(a,&u,&g,&e,&f).\end{array}\end{align}
The isometry $\sigma^*$ can be extended to an action of $(A_2(-1)\oplus U\oplus E_6\oplus E_6\oplus E_6)\otimes \mathbb{Q}$ and so to an action of the overlattice of index $3^2$ of $A_2(-1)\oplus U\oplus E_6\oplus E_6\oplus E_6$ isometric to $\Lambda_{K3}$. In particular $\sigma^*(x)=x$ and $\sigma^*(y)=w+v^{(2)}-v^{(3)}$, with the notation of Theorem \ref{prop: NS(S), H2(S), action of sigma}.

\subsection{The invariant lattice $\Lambda_{K3}^{\sigma^*}$ and its orthogonal complement}
The action of $\sigma^*$ splits $\Lambda_{K3}$ in two sublattices: the invariant sublattice $\Lambda_{K3}^{\sigma^*}$ and its orthogonal complement, denoted by $\Omega_{\Z/3\Z}$ in \cite{GS}. 
By the description of $\sigma^*$, one obtains that the sublattice of $A_2(-1)\oplus U\oplus E_6\oplus E_6\oplus E_6$ invariant for $\sigma^*$ is spanned by the classes $(a_h,0,0,0,0)$, $(0,u_k,0,0,0)$, $(0,0,e_i,e_i,e_i)$ where $h,k=1,2$, $i=1,\ldots, 6$ and $a_h$, $u_k$ and $e_i$ are generators of the lattices  $A_2(-1)$, $U$ and $E_6$ respectively. So  $$\left(A_2(-1)\oplus U\oplus E_6\oplus E_6\oplus E_6\right)^{\sigma^*}\simeq A_2(-1)\oplus U\oplus E_6(3),$$ where the generators of $E_6(3)$ are  $e_i^{(1)}+e_i^{(2)}+e_i^{(3)}$, $i=1,\ldots, 6$. Moreover also the class $x=(\sum_{j=1}^3(e_1^{(j)}+2e_2^{(j)}+e_4^{(j)}+2e_5^{(j)}))/3$ is contained in $(\Lambda_{K3})^{\sigma^*}$. Hence $(\Lambda_{K3})^{\sigma^*}$ is an overlattice of index 3 of $A_2(-1)\oplus U\oplus E_6(3)$, isometric to $A_2(-1)\oplus U\oplus E_6^*(3)$.

Let us now consider the orthogonal complement of $(\Lambda_{K3})^{\sigma^*}$.

\begin{proposition}\label{prop: construction K12 as overlattice}
The lattice $\left(\Lambda_{K3}^{\sigma^*}\right)^{\perp}$ is isometric to the lattice $K_{12}$ and it is spanned by $k_i:=e_i^{(1)}-e_i^{(2)}$, $k_{i+6}:=e_i^{(1)}-e_i^{(3)}$, $i=1,\ldots, 6$ and by the class 
\begin{align}\label{eq: n extra classes K12}z:=\left(k_1+k_4+k_7+k_{10}+2(k_2+k_5+k_8+k_{11})\right)/3.\end{align}
\end{proposition}
\proof We already described $\Lambda_{K3}$ as overlattice of $A_2(-1)\oplus U\oplus E_6\oplus E_6\oplus E_6$ and we described the sublattice $\left(\Lambda_{K3}\right)^{\sigma^*}$ as the lattice generated by $(a_i,0,0,0,0)$, $(0,u_i,0,0,0)$, $(0,0,e_j,e_j,e_j)$, $i=1,2$, $j=1,\ldots,6$, and $x$.
So, $\left(\Lambda_{K3}^{\sigma^*}\right)^{\perp}$ is spanned, at least over $\Q$, by the classes $k_i$, $i=1,\ldots, 12$. The intersection matrix of these classes is the block matrix \begin{eqnarray}\label{eq: K12tilde}\left[\begin{array}{c|c} E_6(2)&E_6\\ \hline E_6&E_6(2)\end{array}\right]\end{eqnarray}

We denote by $\widetilde{K_{12}}$ the lattice whose bilinear form is given by the previous matrix.
Since the determinant of this matrix is $3^8$ and the discriminant of $K_{12}$ is $3^6$, we deduce that $K_{12}$ is an overlattice of index 3 of the lattice $\widetilde{K_{12}}$.
This overlattice is obtained by adding to $\{ k_i\}_{i=1,
\ldots,12}$ another generator. In order to conclude, one has to show that this other generator can be chosen to be $z$. In particular one has to show that $z\in\left(\Lambda_{K3}^{\sigma^*}\right)^{\perp}$. To do this, it suffices to go back to the previous section and to write the  class $z$ in terms of generators of $\NS(S)$, where $S$ is the specific K3 surface considered in Section \ref{sec: a particular surface} (recalling that the three copies of $E_6$ in terms of curves on $S$ are given in \eqref{eq: 3E6}).
The class $z$ corresponds to the class $$2F-2C_1^{(1)}-4C_2^{(1)}-3C_3^{(1)}-2C_4^{(1)}-3C_5^{(1)}-2C_6^{(1)}-2C_2^{(2)}-2C_3^{(2)}-C_4^{(2)}-2C_5^{(2)}-C_6^{(2)},$$
which is contained in $\left(\NS(S)^{\sigma^*}\right)^{\perp}$ (since it is orthogonal to the generators of $\NS(S)^{\sigma^*}$). As a consequence, the lattice $K_{12}\simeq \left(H^2(S,\mathbb{Z})^{\sigma^*}\right)^{\perp}$ is the overlattice of index 3 of $\widetilde{K_{12}}$ obtained by adding $z$ to $\widetilde{K_{12}}$.\endproof

\subsection{The map $\pi_*$}\label{subsec: the map pi_*}
Given a K3 surface $X$ with a symplectic automorphism $\sigma$, the quotient $X/\langle\sigma\rangle$ is a singular surface, whose desingularization is a K3 surface $Y$. Hence there is a generically 3:1 rational map $\pi:X\dashrightarrow Y$, induced by the quotient map $q:X\ra X/\langle\sigma\rangle$.

Hence $\pi_*$ maps the second cohomology group of $X$ to a sublattice of the second cohomology group of $Y$. Since both these surfaces are K3 surfaces, $\pi_*$ maps $\Lambda_{K3}\simeq H^2(X,\Z)$ to a sublattice of $\Lambda_{K3}\simeq H^2(Y,\Z)$. 

\begin{proposition}\label{prop: pi_*} The map $\pi_*$ acts on $A_2(-1)\oplus U\oplus E_6\oplus E_6\oplus E_6$ as follows:
	$$\begin{array}{ccccccccccccccc}\pi_*:&A_2(-1)\oplus &U\oplus &E_6\oplus &E_6\oplus& E_6&\ra& A_2(-3)\oplus &U(3)\oplus &E_6\\
		&(a,&u,&e,&f,&g)&\mapsto& (a,&u,& e+f+g).\end{array}$$
	Its extension to $\Lambda_{K3}$ is such that	
	$\pi_*(H^2(X,\mathbb{Z}))\simeq A_2^*(-3)\oplus U(3)\oplus E_6\simeq A_2(-1)\oplus U(3)\oplus E_6$, which is an overlattice of index 3 of $A_2(-3)\oplus U(3)\oplus E_6$. It is a lattice of rank $10$, signature $(3,7)$ and discriminant group $(\Z/3\Z)^4$.
\end{proposition}
\proof
If $\alpha\in H^2(X,\mathbb{Z})^{\sigma^*}$ is a $\sigma^*$-invariant class, applying the push-pull formula (cf. \cite[Theorem 1.23]{EH}) one obtains $$(\pi_*(\alpha) \cdot \pi_*(\alpha))=3\alpha^2$$
indeed $\pi^*\pi_*(\alpha)=\alpha+\sigma^*(\alpha)+(\sigma^*)^2(\alpha)=3\alpha$.

Given $\alpha_1, \alpha_2\in H^2(X,\mathbb{Z})^{\sigma^*}$ one obtains  $$(\pi_*\alpha_1\cdot\pi_*\alpha_2)=\frac{1}{3}(\pi^*\pi_*\alpha_1\cdot \pi^*\pi_*\alpha_2)=\frac{1}{3}(3\alpha_1\cdot3\alpha_2)=3(\alpha_1\cdot\alpha_2).$$

Hence if $\alpha\in A_2(-1)\oplus U\subset H^2(X,\mathbb{Z})^{\sigma^*}$, then $\pi_*$ multiplies the form by 3, so $\pi_*(A_2(-1)\oplus U)=A_2(-3)\oplus U(3)$.

Let us now consider the image of the classes of the form $(0,0,e,0,0)\in A_2(-1)\oplus U\oplus E_6\oplus E_6\oplus E_6$: let us denote, as before, by $e_i^{(j)}$, $i=1,\ldots, 6$ the basis of the $j$-th copy of $E_6$.
Then
$$(\pi_*e_i^{(1)}\cdot \pi_*e_j^{(1)})=\frac{1}{3}(\pi^*\pi_*e_i^{(1)}\cdot \pi^*\pi_*e_j^{(1)})=\frac{1}{3}((e_i^{(1)}+e_i^{(2)}+e_i^{(3)})\cdot (e_j^{(1)}+e_j^{(2)}+e_j^{(3)}))=(e_i^{(1)}\cdot e_j^{(1)}),$$
where we used: $$\begin{array}{c}\pi^*\pi_*(e_i^{(1)})=e_i^{(1)}+\sigma^*(e_i^{(1)})+(\sigma^*)^2(e_i^{(1)})=e_i^{(1)}+e_i^{(2)}+e_i^{(3)},\\ (e_i^{(h)}\cdot e_j^{(h)})=(e_i^{(1)}\cdot e_j^{(1)}),\mbox{ and } (e_i^{(h)}\cdot e_j^{(k)})=0,\mbox{ if } h\neq k.\end{array}$$

Hence $\pi_*(E_6\oplus E_6\oplus E_6)\simeq E_6.$

So we obtain $\pi_*(A_2(-1)\oplus U\oplus E_6\oplus E_6\oplus E_6)\simeq A_2(-3)\oplus U(3)\oplus E_6.$ 
In order to find $\pi_*(H^2(X,\mathbb{Z}))$, it remains to understand the images of the classes $x$ and $y$.  

The class $x=\left(\sum_j(e_1^{(j)}+2e_2^{(j)}+e_4^{(j)}+2e_5^{(j)})\right)/3$ is mapped to $$\pi_*(x)=\pi_*(e_1^{(1)}+2e_2^{(1)}+e_4^{(1)}+2e_5^{(1)})$$
since $\pi_*(e_i^{(1)}+e_i^{(2)}+e_i^{(3)})=3\pi_*(e_i^{(1)})$.

The class $y=\left(a_1+2a_2+e_1^{(1)}+2e_2^{(1)}+e_4^{(1)}+2e_5^{(1)}-e_1^{(2)}-2e_2^{(2)}-e_4^{(2)}-2e_5^{(2)}\right)/3$ is mapped to $$\pi_*(y)=\pi_*((a_1+2a_2)/3)=(\pi_*a_1+2\pi_*a_2)/3$$ since $\pi_*e_i^{(1)}-\pi_*e_i^{(2)}=\pi_*e_i^{(1)}-\pi_*e_i^{(1)}=0$.
The vectors $\pi_*a_i$ are the generators of $A_2(-3)$ and
hence we are constructing an overlattice of index 3 of $A_2(-3)$, isometric to $A_2^*(-3)$. This lattice is isometric to $A_2(-1)$, with basis
$$a_1':=(\pi_*a_1+2\pi_*a_2)/3,\ \  a_2':=\pi_*(a_2)-\left((\pi_*a_1+2\pi_*a_2)/3\right)=(-\pi_*a_1+\pi_*a_2)/3.$$
We observe that $a_1'=\pi_*(y)$ and $a_2'=\pi_*(a_2)-\pi_*(y)$.


\endproof

\subsection{The cohomology of the quotient K3 surface $Y$}\label{subsec: Cohomology of Y and M3}
Let $X$ be a K3 surface admitting a symplectic automorphism $\sigma$ of order 3, then the surface $X/\langle\sigma\rangle$ 
has 6 singularities of type $A_2$ (see \cite[Section 5]{Nik}). We denote by $Y$ the desingularization of $X/\langle\sigma\rangle$, which introduces 12 irreducible curves, 6 disjoint pairs of rational curves meeting in a point. We call the resolution map $\beta:Y\ra X/\langle\sigma\rangle$ and we denote by $M_i^{(j)}$, $i=1,2$, $j=1,\ldots, 6$ the curves introduced by $\beta$, with the following conventions: $$M_1^{(j)}M_2^{(j)}=1,\ \ \left(M_i^{(j)}\right)^2=-2\mbox{ and } M_i^{(j)}M_{h}^{(k)}=0\mbox{ if }j\neq k.$$
By \cite[Section 6]{Nik} the minimal primitive sublattice of $H^2(Y,\Z)$ which contains the curves $M_i^{(j)}$ contains also the class $\hat{M}:=\sum_{j=1}^6(M_1^{(j)}+2M_2^{(j)})/3$. Thus it is an overlattice of index 3 of $A_2^{\oplus 6}$ and so it is a negative definite lattice of rank 12 and discriminant group $(\Z/3\Z)^4$. This lattice is denoted by $M_{\Z/3\Z}$.

The lattice $\pi_*(H^2(X,\mathbb{Z}))$ is naturally embedded in $H^2(Y,\mathbb{Z})$. The curves introduced by $\beta$ are orthogonal to $\beta^*(D)$ for each divisor $D\in \pi_*(\NS(X))$ and of course they are also orthogonal to each class in $\pi_*(T_X)$. So we have the following orthogonal decomposition (which holds over the rational field, but not over the integer ring) 
\begin{equation}\label{eq: H^2Y direct sum pi* and M3}H^2(Y,\mathbb{Q})=\left(M_{\Z/3\Z}\oplus \pi_*(H^2(X,\Z))\right)\otimes \Q.\end{equation}
\subsection{Gluing $\pi_*(H^2(X,\Z))$ and $M_{\Z/3\Z}$ to obtain $H^2(Y,\Z)$.}

By \eqref{eq: H^2Y direct sum pi* and M3}, the orthogonal complement of $M_{\Z/3\Z}$ in $H^2(Y,\mathbb{Z})$ is either $\pi_*(H^2(X,\mathbb{Z}))$ or an overlattice of finite index of $\pi_*(H^2(X,\mathbb{Z}))$.

By definition, $M_{\Z/3\Z}$, is primitively embedded in $H^2(Y,\mathbb{Z})$.  Since $H^2(Y,\Z)$ is unimodular, the discriminant of the orthogonal complement of $M_{\Z/3\Z}$ in $H^2(Y,\Z)$ is $-d(M_{\Z/3\Z})$. So, if $d(\pi_*((H^2(X,\mathbb{Z}))))=-d(M_{\Z/3\Z})$, then $\pi_*((H^2(X,\mathbb{Z})))$ is the orthogonal complement of $M_{\Z/3\Z}$. One immediately checks that $|d(M_{\Z/3\Z})|=3^4=|d(\pi_*(H^2(X,\mathbb{Z})))|$, so
$\pi_*((H^2(X,\mathbb{Z})))$ is the orthogonal complement of $M_{\Z/3\Z}$.

 We recall that the orthogonal complement of a given sublattice inside a bigger lattice $\Lambda$, is primitively embedded in $\Lambda$. In particular, $\pi_*(H^2(Y,\Z))\simeq M_{\Z/3\Z}^{\perp}$ is primitively embedded in $H^2(Y,\Z)$. 
 Therefore to construct the unimodular lattice $H^2(Y,\Z)$ one has to glue these two primitive sublattices $\pi_*(H^2(Y,\Z))$ and $M_{\Z/3\Z}$. To this purpose we need the following description of the discriminant group of $M_{\Z/3\Z}$.


\begin{lemma}\label{lemma: base discriminant M3} With the notation of Section \ref{subsec: Cohomology of Y and M3}, the discriminant group of $M_{\Z/3\Z}$ is generated by $$\begin{array}{ll}b_1:=z_1+z_2+z_3,& b_2:=z_1+z_2+z_4,\\b_3:=z_2-z_3+z_4-z_5,& b_4:=-z_1+z_3-z_4+z_5,\end{array}$$ where $$z_j:=\left(M_1^{(j)}+2M_2^{(j)}\right)/3.$$\end{lemma}
\proof It suffices to check that $b_kM_i^{(j)}\in\Z$ and that $b_i$ are independent in $(\Z/3\Z)^4\simeq A_{M_{\Z/3\Z}}$. The latter statement can be checked by observing that the discriminant form computed on $\{b_1,b_2,b_3,b_4\}$ is non degenerate, and indeed it is the opposite of the discriminant form of $U(3)\oplus A_2(-1)\oplus E_6$.\endproof 

\begin{proposition}\label{prop: H^2(Y) extra classes} Denoted by $a_i', u_i', e_j'$, $i=1,2$, $j=1,\ldots, 6$ the standard generators of $A_2(-1)$, $U(3)$ and $E_6$ respectively, the overlattice $H^2(Y,\mathbb{Z})$ of $\pi_*(H^2(X,\Z))\oplus M_{Z/3\Z}\simeq A_2(-1)\oplus U(3)\oplus E_6\oplus M_{\Z/3\Z}$ is obtained by adding the classes $$\begin{array}{llll}n_1:=&\frac{a_1'+2a_2'}{3}+b_3
&n_2:=&\frac{e_1'+2e_2'+e_4'+2e_5'}{3}+b_4\\
n_3:=&\frac{u_1'}{3}+b_1&n_4:=&\frac{u_2'}{3}+b_2\end{array}$$
to $A_2(-1)\oplus U(3)\oplus E_6\oplus M_{\Z/3\Z}$.
\end{proposition}
\proof It suffices to check that the lattice obtained by adding to the generators of $A_2(-1)\oplus U(3)\oplus E_6\oplus M_{\Z/3\Z}$ the classes $n_1$, $n_2$, $n_3$ and $n_4$ is even and unimodular. It follows that it is an even unimodular lattice of signature $(3,19)$ and then it is isometric to the lattice $\Lambda_{K3}\simeq H^2(Y,\Z)$.\endproof

Here we are interested in describing the classes $n_i$ in terms of the curves contained $M_{\Z/3\Z}$. In order to construct the overlattice of $A_2(-1)\oplus U(3)\oplus E_6\oplus M_{\Z/3\Z}$, it is equivalent if one adds to it the class $n_i$ or the class $n_i+\sum_{i}\alpha_{i,j}M_i^{(j)}$, $\alpha_{i,j}\in\mathbb{Z}$, indeed these two classes are equivalent in the discriminant group. So we can rewrite the classes $n_i$ in the following form 
$$n_1\sim\frac{2\pi_*(a_2)-\pi_*(y)+M_1^{(2)}+2M_2^{(2)}+2M_1^{(3)}+M_2^{(3)}+M_1^{(4)}+2M_2^{(4)}+2M_1^{(5)}+M_2^{(5)}}{3},$$
$$n_2\sim\frac{\pi_*(e_1+2e_2+e_4+2e_5)+2M_1^{(1)}+M_2^{(1)}+M_1^{(3)}+2M_2^{(3)}+2M_1^{(4)}+M_2^{(4)}+M_1^{(5)}+2M_2^{(5)}}{3},$$
$$n_3\sim \frac{\pi_*(u_1)+M_1^{(1)}+2M_2^{(1)}+M_1^{(2)}+2M_2^{(2)}+M_1^{(3)}+2M_2^{(3)}}{3},$$
$$ n_4\sim \frac{\pi_*(u_2)+M_1^{(1)}+2M_2^{(1)}+M_1^{(2)}+2M_2^{(2)}+M_1^{(4)}+2M_2^{(4)}}{3}, $$
where we used that $a_1'=\pi_*(y)$ and $a_2'=\pi_*(a_2)-\pi_*(y)$, by proof of Proposition \ref{prop: pi*}.
\subsection{The map $\pi^*$}\label{subsec: pi*} We consider the map $\pi^*$, dual of the map $\pi_*$ considered above.
\begin{proposition}\label{prop: pi*}
	The map $\pi^*$ acts as follows
	$$\begin{array}{cccccccccc}\pi^*:&A_2(-3)\oplus &U(3)\oplus& E_6&\ra& A_2(-1)\oplus &U\oplus &E_6\oplus& E_6\oplus& E_6\\
		&(\alpha,&\mu,&e)&\mapsto& (3\alpha,&3\mu,&e,&e,&e).\end{array}$$
	With the notation of Proposition \ref{prop: H^2(Y) extra classes},  $H^2(Y,\Z)\simeq \Lambda_{K3}$ is obtained by adding to $A_2(-3)\oplus U(3)\oplus E_6\oplus M_{\Z/3\Z}$ the classes $a_1'$ and $n_h$, $h=1,2,3,4$, so the extension of $\pi^*$ to $H^2(Y,\Z)$ is given by
	$$\pi^*(a_1')=a_1+2a_2,\ \ \pi^*(n_1)=2a_2-y,\ \ \pi^*(n_2)=x,\ \ \pi^*(n_3)=u_1,\ \ \pi^*(n_4)=u_2$$
	where $a_i$ and $u_i$ are generators of $A_2(-1)$ and $U$ respectively.	 
\end{proposition}
\proof 
 By Proposition \ref{prop: pi_*}, $\pi_*(\gamma)=\gamma$ if $\gamma\in A_2(-1)\oplus U\subset H^2(X,\Z)^{\sigma^*}$ and the map  ${\pi_*}_{|H^2(X,\Z)^{\sigma^*}}$ multiplies the form by 3. For every $\gamma\in A_2(-1)\oplus U\subset H^2(X,\mathbb{Z})^{\sigma^*}$ and $\alpha\in \pi_*(H^2(X,\Z)^{\sigma^*})\simeq A_2(-3)\oplus U(3)\oplus E_6$, it holds $(\pi^*\alpha\cdot\gamma)=(\alpha\cdot\pi_*\gamma)$. So we obtain 
 $$(\pi^*\alpha\cdot\gamma)=(\alpha\cdot\pi_*\gamma)=3(\alpha\cdot \gamma)$$
and hence $\pi^*(\alpha)=3\alpha$ for every $\alpha\in A_2(-3)\oplus U(3)$. 

Let $e\in E_6\subset H^2(Y,\Z)$ and $(f,0,0)\in E_6^{\oplus 3}\subset H^2(X,\Z)$. By push-pull formula and Proposition \ref{prop: pi_*}, $$(\pi^*e\cdot (f,0,0))=(e\cdot \pi_*(f,0,0))=(e\cdot f).$$
Similarly $(\pi^*e\cdot (0,f,0))=(e\cdot f)$ and $(\pi^*e\cdot (0,0,f))=(e\cdot f)$.
Hence $$\pi^*(e)=(e,e,e)\in E_6\oplus E_6\oplus E_6\subset H^2(X,\Z).$$ 

Now we want to extend this map to the lattice $A_2(-1)\oplus U(3)\oplus E_6$, obtained by adding to $A_2(-3)\oplus U(3)\oplus E_6$ the class $a_1' (=\pi_*((a_1+2a_2)/3)=\pi_*(y))$. Since $\pi^*(\alpha)=3\alpha$ for every $\alpha\in A_2(-3)$,  we obtain $\pi^*(a_1')=\pi^*(\pi_*(a_1+2a_2)/3)=a_1+2a_2$.

In order to extend this map to $H^2(Y,\Z)\simeq \Lambda_{K3}$, 
it remains to compute $\pi^*(n_h)$ for $h=1,2,3,4$. We observe that $\pi^*(M_i^{(j)})=0$ indeed if one reconstruct $X$ from $Y$, one considers the (minimal resolution of the) triple cover $Y$ branched on the curves $M_i^{(j)}$. This surface is not minimal and then one contracts the $(-1)$-curves. Among these contracted curves, one finds the curves which are covers of the curves $M_i^{(j)}$. In particular there are no curves on $X$ which are the inverse image of curves $M_i^{(j)}\subset Y$. So 
$$\begin{array}{l}\pi^*(n_1)=\pi^*\left(2\pi_*(a_2)-\pi_*(y)\right)/3=(2\pi^*\pi_*(a_2)-\pi^*\pi_*(y))/3=2a_2-y\\
\pi^*(n_2)=\sum_{j=1}^3\left(e_1^{(j)}+2e_2^{(j)}+e_4^{(j)}+2e_5^{(j)}\right)/3=x,\ \ 
\pi^*(n_3)=u_1,\ \ 
\pi^*(n_4)=u_2.\\
\end{array}$$
\endproof

\section{Families of projective K3 surfaces admitting a symplectic automorphism of order 3 and of their quotients}\label{sec: projective families}
The purpose of this Section is to describe the families of the projective K3 surfaces admitting a symplectic automorphism of order 3 and the families of their quotients in terms of families of lattices polarized K3 surfaces. So, we want to reduce the geometric property ``to admit a symplectic automorphism of order 3" (resp. ``to be a quotient of a K3 surface by a symplectic automorphism of order 3") into a lattice theoretic property. It is already known that this is possible by considering overlattices of the lattices $K_{12}$ and $M_{\Z/3\Z}$ with certain properties (see \cite{Nik},  \cite{GS}, \cite{G}). Here we classify these overlattices providing a complete description of the two families considered, see Theorems \ref{theor:NS X} and \ref{theor:NS Y} and Corollaries  \ref{cor: families X} and \ref{cor: families Y}.

\subsection{Projective K3 surfaces with a symplectic automorphism of order 3}

\begin{lemma}\label{lemm: orbits AK12} Let $q_{A_{K_{12}}}$ be the quadratic form on the discriminant group $A_{K_{12}}\simeq\left(\Z/3\Z\right)^6$, then $q_{A_{K_{12}}}(v)$ for $v\in A_{K_{12}}$ is one of the following three values $0,\frac{2}{3}, \frac{4}{3}$.  
	
	The action induced by $O(K_{12})$ on $A_{K_{12}}$ has four orbits: the orbit which contains 0, and the orbits:\begin{itemize}\item
		$o_0=\{v\in A_{K_{12}}\mbox{ such that }q_{A_{K_{12}}}(v)=0\mbox{ and }v\neq 0\}$;
		\item 
		$o_1=\{v\in A_{K_{12}}\mbox{ such that }q_{A_{K_{12}}}(v)=\frac{2}{3}\}$;
		\item$o_2=\{v\in A_{K_{12}}\mbox{ such that }q_{A_{K_{12}}}(v)=\frac{4}{3}\}$.
	\end{itemize}
\end{lemma}
\proof 
The result follows by  \cite[Section 3, in particular formula (10)]{CS}. Nevertheless we give an idea of the proof. The isometries of the $K_{12}$ lattice induce isometries of its discriminant group.
If two vectors $a\neq 0$ and $b\neq 0$ in $A_{K_{12}}$  are in the same orbit for the action of $O(A_{K_{12}})$, then $q_{A_{K_{12}}}(a)=q_{A_{K_{12}}}(b)$. In order to reverse the implication, i.e. to prove that the orbit of a non-zero vector is completely determined by the value of the discriminant form on it, it suffices to find an isometry in $O(A_{K_{12}})$ which maps $a$ in $b$ if $q_{A_{K_{12}}}(a)=q_{A_{K_{12}}}(b)$ and which is induced by an isometry of $O(K_{12})$. 
The vectors $$\begin{array}{c}g_1:=(k_7+2k_8+k_{10}+2k_{11})/3,\ g_2:=\left(k_6+k_{12}\right)/3,\\ g_3:=\left(k_5+k_{11}\right)/3,\ g_4:=\left(k_4+k_{10}\right)/3,\ g_5:=\left(k_3+k_{9}\right)/3,\ g_6:=\left(k_2+k_{8}\right)/3\end{array}$$ 
form a basis of $A_{K_{12}}$ such that the quadratic form $q_{A_{K_{12}}}$ is 
$$q_{A_{K_{12}}}(v)=\frac{4}{3}x_1^2 + \frac{2}{3}(x_2^2+x_3^2+x_4^2+x_5^2+x_6^2)+\frac{4}{3}(x_2x_5 + x_3x_4+x_4x_5+x_5x_6),$$
where $v=x_1g_1 +\ldots+x_6g_6$. Let us consider the following order 2 isometry of $K_{12}$: $k_1\stackrel{\varphi}{\leftrightarrow} k_5$, $k_2\stackrel{\varphi}{\leftrightarrow} k_4$, $k_7\stackrel{\varphi}{\leftrightarrow}k_{11}$, $k_8\stackrel{\varphi}{\leftrightarrow} k_{10}$, $\varphi(k_i)=k_i$, $i=3,6,9,12$.
Then $\varphi\in O(K_{12})$ induces $\varphi_A\in O(A_{K_{12}})$ which  acts as follows 
$$\varphi_A(v)=-x_1g_1+x_2g_2+x_3g_3+(-x_3+x_6)g_4+x_5g_5+(x_3+x_4)g_6.$$
By considering $\langle\varphi_A, -id\rangle\in O(A_{K_{12}})$ one shows that the vectors $g_3, -g_3, g_3-g_4+g_6, -g_3+g_4-g_6$ are in the same orbit for $O(A_{K_{12}})$. 
A complete table of all the isometries needed to show the result can be found in \cite[Appendix 1]{P} and is obtained by using the integral lattices package of Sage, \cite{Sage}.
\endproof


\begin{proposition}\label{prop:overlattice 2d+K12}
	Let $d$ be a positive integer. 
	
	If $d\not\equiv 0\mod 3$, then there exist no even overlattices of finite index of  $\langle 2d\rangle\oplus K_{12}$ such that $\langle 2d\rangle$ and $K_{12}$ are primitively embedded in it.
	
	If $d\equiv 0\mod 3$, then there exists a unique overlattice of finite index of $\langle 2d\rangle\oplus K_{12}$ such that $\langle 2d\rangle$ and $K_{12}$ are primitively embedded in it. The index is 3 and the generators of the overlattice are the ones of $\langle 2d\rangle\oplus K_{12}$ and $\frac{L+g}{3}$ where $L$ is a generator of $\langle 2d\rangle$ and $g\in K_{12}$ can be chosen as follows:

if $d\equiv 0\mod 9$, then $g=k_1+k_3+k_5+k_7+k_9+k_{11}$

		 if $d\equiv 3\mod 9$, then $g=k_1+k_3+k_7+k_9;$

		 if $d\equiv 6\mod 9$, then $g=k_1+k_7.$
	 
\end{proposition}
\proof Let $\Lambda$ be a proper overlattice of finite index of $\langle 2d\rangle \oplus K_{12}$ in which $K_{12}$ and $\langle2d\rangle$ are primitively embedded. Let $f$ be one generator of $\Lambda/(\langle 2d\rangle\oplus K_{12})$ then $f\in A_{(\langle 2d\rangle\oplus K_{12})}$ and hence $$f=(\alpha L+\beta g)/t,\mbox{ with }\alpha L/t\in A_{\langle 2d\rangle},\ \beta g/t\in A_{K_{12}}.$$
Since $A_{K_{12}}\simeq (\Z/3\Z)^6$, $t$ has to be equal to 3 and hence $\alpha,\beta\in \{0,1,2\}$. 

Since $\langle 2d\rangle=\Z L$ (resp. $K_{12}$) is primitively embedded in $\Lambda$, $\beta\neq 0$ (resp. $\alpha\neq 0$). 
Moreover, since $\Lambda/(\langle 2d\rangle\oplus K_{12})\subset (\Z/3\Z)^c$ one can freely substitute $g$ (resp. $f$) by $2g$ (resp. $2f$) because 2 is invertible in $\Z/3\Z$. Indeed an easy computation shows that one can assume that $\beta=1$ (resp. $\alpha=1$) and $g/3\in A_{K_{12}}$.

So $f=(L+g)/3$. The bilinear form on $\Lambda$ takes integer values and the lattice is even, which implies
$$fL=2d/3\in\Z\Leftrightarrow d\equiv 0\mod 3,\ \ f^2=(2d+g^2)/9\in 2\Z\Leftrightarrow 2d\equiv -g^2\mod 18.$$

In particular the existence of $f$ implies that $d\equiv 0\mod 3$. Moreover 

if $d\equiv 3\mod 9$, then $g/3\in A_{K_{12}}$ satisfies $q_{A_{K_{12}}}(\frac{g}{3})=-\frac{2}{3},$

if $d\equiv 6\mod 9$, then $g/3\in A_{K_{12}}$ satisfies $q_{A_{K_{12}}}(\frac{g}{3})=-\frac{4}{3},$

if $d\equiv 0\mod 9$, then $g/3\in A_{K_{12}}$ satisfies $q_{A_{K_{12}}}(\frac{g}{3})=0.$

By Lemma \ref{lemm: orbits AK12} one can choose arbitrarily an element $g\in K_{12}$ such that $g/3\in K_{12}^{\vee}/K_{12}$ is non trivial and $q_{A_{12}}(g)$ is the required one, since the choice of the value of $q_{A_{K_{12}}}(g)$ determines a unique orbit. In the statement we made a specific choice for $g$, any other is equivalent up to isometries of $K_{12}$.

It remains to prove that $\Lambda/(\langle 2d\rangle\oplus K_{12})$ is a cyclic group and thus the index of $\langle 2d\rangle\oplus K_{12}\hookrightarrow \Lambda$ is 3. Suppose there exist two independent generators $f_1$ and $f_2$ of $\Lambda/(\langle 2d\rangle\oplus K_{12})$. Both of them can be chosen to be of the form $f_i=(L+g_i)/3$. So $f_1-f_2=(g_1-g_2)/3\in \Lambda$ and $(g_1-g_2)/3\in A_{K_{12}}$ is non trivial. Then $K_{12}$ is not primitively embedded in $\Lambda$, which is absurd.\endproof

\begin{definition}\label{def: overlattice}
	If $d\equiv 0 \mod 3$ we denote by $\left(\langle 2d\rangle\oplus K_{12}\right)'$ the overlattice constructed in the previous proposition.
\end{definition}

\begin{theorem}{\rm(See also \cite[Proposition 5.1]{GS})}\label{theor:NS X}
	Let $X$ be a K3 surface. Then $X$ admits a symplectic automorphism of order 3 if and only if $K_{12}$ is primitively embedded in $\NS(X)$. If $X$ is a projective K3 surface which admits a symplectic automorphism of order 3, then $\rho(X)\geq 13$. 
	Assume $\rho(X)=13$ and let $L$ be a generator of $K_{12}^{\perp_{NS(X)}}$, so $L^2=2d$, $d\in\mathbb{N}_{>0}$. Then\begin{itemize}
		\item if $d\not\equiv 0\mod 3$, then $\NS(X)=\langle 2d\rangle\oplus K_{12}$;
		\item if $d \equiv 0\mod 3$, then $\NS(X)$ is either $\langle 2d\rangle\oplus K_{12}$ or $(\langle 2d\rangle\oplus K_{12})'$.
	\end{itemize} 
\end{theorem}
\proof The statement was already proved in \cite[Proposition 5.1]{GS}, the unique difference is that in \cite[Proposition 5.1]{GS}, one proved that if $d\equiv 0\mod 3$ then $\NS(X)$ could be an overlattice of index 3 of $\langle 2d\rangle\oplus K_{12}$. By Proposition \ref{prop:overlattice 2d+K12}, we know that this overlattice exists, it is unique up to isometries and it coincides with the lattice of the Definition \ref{def: overlattice}. 
\endproof

\begin{corollary}{\rm(See \cite[Proposition 5.2]{GS})}\label{cor: families X}
	Let $X$ be a $\left(\langle 2d\rangle\oplus K_{12}\right)$-polarized K3 surface or $\left(\langle 2d\rangle\oplus K_{12}\right)'$-polarized K3 surface. Then $X$ is projective and admits a symplectic automorphism of order 3. Let $$\mathcal{S}=\{\mbox{projective K3 surfaces admitting an order 3 symplectic automorphism}\}/\simeq$$ where $\simeq$ is the isomorphisms of polarized K3 surfaces. Then $\mathcal{S}$ is  
	$$\bigcup_{d\in\mathbb{N}_{>0}}\left(\{\left(\langle 2d\rangle\oplus K_{12}\right)-\mbox{ polarized K3s }\}_{/\simeq}\right)\bigcup_{{\footnotesize{\begin{array}{c}d\in\mathbb{N}_{>0},\\d\equiv 0\mod 3\end{array}}}}\left(\{\left(\langle 2d\rangle\oplus K_{12}\right)'-\mbox{ polarized K3s} \}_{/\simeq}\right).$$
	In particular $\mathcal{S}$ is the union of countably many components all of dimension 7.
\end{corollary}

\subsection{Projective K3 surfaces quotients of K3 surfaces by a symplectic automorphism of order 3}

\begin{lemma}\label{lemm: orbits AM3} Let $q_{A_{M_{\Z/3\Z}}}$ be the quadratic form of the discriminant group $A_{M_{\Z/3\Z}}\simeq\left(\Z/3\Z\right)^4$, then for $v\in A_{M_{\Z/3\Z}}$, $q_{A_{M_{\Z/3\Z}}}(v)$ is one of the following three values $0,\frac{2}{3}, \frac{4}{3}$.  
	
	The actions induced by $O(M_{\Z/3\Z})$ on $A_{M_{\Z/3\Z}}$ has four orbits: the orbit which contains 0, and the orbits:\begin{itemize}\item
		$p_0=\{v\in A_{M_{\Z/3\Z}}\mbox{ such that }q_{A_{M_{\Z/3\Z}}}(v)=0\mbox{ and }v\neq 0\}$;
		\item 
		$p_1=\{v\in A_{M_{\Z/3\Z}}\mbox{ such that }q_{A_{M_{\Z/3\Z}}}(v)=\frac{2}{3}\}$;
		\item $p_2=\{v\in A_{M_{\Z/3\Z}}\mbox{ such that }q_{M_{\Z/3\Z}}(v)=\frac{4}{3}\}$.
	\end{itemize}
\end{lemma}
\proof We sketch the proof which is analogous to the one of Lemma \ref{lemm: orbits AK12}. The quadratic form $q_{A_{M_{\Z/3\Z}}}$ is
$q_{A_{M_{\Z/3\Z}}}(v)=-\frac{2}{3}(x_3^2+x_4^2)-\frac{1}{3}(2x_1x_2+2x_2x_3-2x_2x_4-x_3x_4)$ where $v=x_1b_1+x_2b_2+x_3b_3+x_4b_4$ and $b_i$ are the generators of $A_{M_{\Z/3\Z}}$ defined in Lemma \ref{lemma: base discriminant M3}. As example, we observe that the isometry of $M_{\Z/3\Z}$ which permutes $M_i^{(3)}$ with $M_i^{(4)}$, $i=1,2$, and acts as the identity on $M_i^{(j)}$ for $i=1,2$, $j=1,2,5,6$, induces on $A_{M_{\Z/3\Z}}$ an isometry which maps $b_1$ to $b_2$.
	%
\endproof

\begin{proposition}\label{prop:overlattice 2d+M3}{\rm (See \cite[Proposition 5.5]{G})}
	Let $e$ be a positive integer. 
	
	If $e\not\equiv 0\mod 3$, then there exists no even overlattices of finite index of  $\langle 2e\rangle\oplus M_{\Z/3\Z}$ such that $\langle 2e\rangle$ and $M_{\Z/3\Z}$ are primitively embedded in it.
	
	If $e\equiv 0\mod 3$, then there exists a unique overlattice of finite index of $\langle 2e\rangle\oplus M_{\Z/3\Z}$ such that $\langle 2e\rangle$ and $M_{\Z/3\Z}$ are primitively embedded in it. The index is 3 and the generators of the overlattice are the ones of $\langle 2e\rangle\oplus M_{\Z/3\Z}$ and 
	$\frac{H+g}{3}$ where $H$ is a generator of $\langle 2e\rangle$ and $g\in M_{\Z/3\Z}$ can be chosen as follows:

	if $e\equiv 0\mod 9$, then $g=\sum_{i=1}^3\left(M_1^{(i)}+2M_2^{(i)}\right)$

	if $e\equiv 3\mod 9$, then $g=\sum_{i=1}^2\left(2M_1^{(i)}+M_2^{(i)}\right)+\sum_{j=3}^4\left(M_1^{(j)}+2M_2^{(j)}\right)$
	
	if $e\equiv 6\mod 9$, then $g=M_1^{(1)}+2M_2^{(1)}+2M_1^{(2)}+M_2^{(2)}.$
\end{proposition}
\proof The proof is analoguos to the one of Proposition \ref{prop:overlattice 2d+M3} and it is based on Lemma \ref{lemm: orbits AM3}.\endproof

\begin{definition}
	If $e\equiv 0 \mod 3$ we denote by $\left(\langle 2e\rangle\oplus M_{\Z/3\Z}\right)'$ the overlattice constructed in the previous proposition.
\end{definition}

\begin{theorem}\label{theor:NS Y}{\rm (See \cite[Proposition 5.5]{G})}
	Let $Y$ be a K3 surface. It is the minimal resolution of the quotient of a K3 surface by a symplectic automorphism of order 3 if and only if $M_{\Z/3\Z}$ is primitively embedded in $\NS(Y)$. If moreover $Y$ is projective, then $\rho(Y)\geq 13$. 
	Assume $\rho(Y)=13$ and let $H$ be a generator of $M_{\Z/3\Z}^{\perp_{NS(Y)}}$, so $H^2=2e$, $e\in\mathbb{N}_{>0}$. Then\begin{itemize}
		\item if $e\not\equiv 0\mod 3$, then $\NS(Y)=\langle 2e\rangle\oplus M_{\Z/3\Z}$;
		\item if $e\equiv 0\mod 3$, then $\NS(Y)$ is either $\langle 2e\rangle\oplus M_{\Z/3\Z}$ or $(\langle 2e\rangle\oplus M_{\Z/3\Z})'$
	\end{itemize} 
\end{theorem}
\proof The proof is analogous to the one of Theorem \ref{theor:NS X} and it is based on Proposition \ref{prop:overlattice 2d+M3}. \endproof

\begin{corollary}\label{cor: families Y}
	Let $Y$ be a $\left(\langle 2d\rangle\oplus M_{\Z/3\Z}\right)$-polarized K3 surface or a $\left(\langle 2d\rangle\oplus M_{\Z/3\Z}\right)'$-polarized K3 surface. Then $Y$ is projective and it is the desingularization of the quotient of a K3 surface by a symplectic automorphism of order 3. Let $$\mathcal{T}=\{\mbox{projective K3s quotient of K3s by an order 3 symplectic automorphism}\}/\simeq$$ where $\simeq$ is the isomorphisms of K3 surfaces. Then $\mathcal{T}$ is  
	$$\bigcup_{e\in\mathbb{N}_{>0}}\left(\{\left(\langle 2e\rangle\oplus M_{\Z/3\Z}\right)-\mbox{ polarized K3s }\}_{/\simeq}\right)\ \bigcup_{{\footnotesize{\begin{array}{c}e\in\mathbb{N}_{>0},\\e\equiv 0\mod 3\end{array}}}}\left(\{\left(\langle 2e\rangle\oplus M_{\Z/3\Z}\right)'-\mbox{ polarized K3s} \}_{/\simeq}\right).$$
	In particular $\mathcal{T}$ is the union of countably many components all of dimension 7.
\end{corollary}

We now introduce some divisors which will be interesting in the following.

\begin{rem}\label{rem: defi Di}{\rm
Let $e\equiv 0\mod 9$ and $Y$ be a K3 surface such that $\NS(Y)\simeq\left(\langle2e\rangle\oplus M_{\Z/3\Z}\right)'$. Then the following divisors are contained in $\NS(Y)$:
\begin{eqnarray}\label{eq: D_i case tilde 1}\begin{array}{c}D_1:=\frac{H-(M_1^{(1)}+2M_2^{(1)}+M_1^{(2)}+2M_2^{(2)}+M_1^{(3)}+2M_2^{(3)})}{3}\\ \\
 
D_2:=\frac{H-(2M_1^{(4)}+M_2^{(4)}+2M_1^{(5)}+M_2^{(5)}+2M_1^{(6)}+M_2^{(6)})}{3}\\ \\
D_3:=\frac{H-\left(\sum_{i=1}^3(2M_1^{(i)}+M_2^{(i)})+\sum_{j=4}^6(M_1^{(j)}+2M_2^{(j)})\right)}{3}.\end{array}\end{eqnarray}
Indeed the divisor $D_1$ is contained in $\NS(Y)$ by Proposition \ref{prop:overlattice 2d+M3}, the divisor $D_2$ (resp. $D_3$) is obtained from $D_1$ by adding first the class $2\hat{M}$ (resp. $\hat{M}$) and then a linear combination with integer coefficients of the curves $M_i^{(j)}$, i.e. by adding to $D_1$ classes in $M_{\Z/3\Z}\subset \NS(Y)$.

Similarly, if
$e\equiv 3\mod 9$ and $Y$ is a K3 surface such that $\NS(Y)\simeq\left(\langle2e\rangle\oplus M_{\Z/3\Z}\right)'$, the following divisors are contained in $\NS(Y)$:
\begin{equation}\label{eq: D_i case tilde 2}
	\begin{split}D_1:=\left(H-\left(2M_1^{(1)}+M_2^{(1)}+2M_1^{(2)}+M_2^{(2)}+M_1^{(3)}+2M_2^{(3)}+M_1^{(4)}+2M_2^{(4)}\right)\right)/3\\
D_2:=\left(H-\left(M_1^{(1)}+2M_2^{(1)}+M_1^{(2)}+2M_2^{(2)}+M_1^{(5)}+2M_2^{(5)}+M_1^{(6)}+2M_2^{(6)}\right)\right)/3\\ 
D_3:=\left(H-\left(2M_1^{(3)}+M_2^{(3)}+2M_1^{(4)}+M_2^{(4)}+M_1^{(5)}+2M_2^{(5)}+M_1^{(6)}+2M_2^{(6)}\right)\right)/3.\end{split}\end{equation}

If
$e\equiv 6\mod 9$ and $Y$ is a K3 surface such that $\NS(Y)\simeq\left(\langle2e\rangle\oplus M_{\Z/3\Z}\right)'$, the following divisors are contained in $\NS(Y)$:
\begin{equation}\label{eq: D_i case tilde 3}\begin{split}
D_1:=\left(H-\left(M_1^{(1)}+2M_2^{(1)}+2M_1^{(2)}+M_2^{(2)}\right)\right)/3\\ D_2:=\left(H-\left(M_1^{(2)}+2M_2^{(2)}+\sum_{i=3}^6(2M_1^{(i)}+M_2^{(i)})\right)\right)/3\\ D_3:=\left(H-\left(2M_1^{(1)}+M_2^{(1)}+\sum_{i=3}^6(M_1^{(i)}+2M_2^{(i)})\right)\right)/3.\end{split}\end{equation}

If $e\not\equiv 0$ and $Y$ is a K3 surface such that $\NS(Y)\simeq\langle2e\rangle\oplus M_{\Z/3\Z}$, the following divisors are contained in $\NS(Y)$:
\begin{equation}\label{eq: D_i not tilde}D_1:=H,\ \  D_2:=H-\left(\sum_{i=i}^6(2M_1^{(i)}+M_2^{(i)})\right)/3,\ \   D_3:=H-\left(\sum_{i=i}^6(M_1^{(i)}+2M_2^{(i)})\right)/3.\end{equation}
}\end{rem}
\section{Relations between the families of projective K3 surfaces admitting a symplectic automorphism and of their quotients}\label{sec: relations between families}
The main results of this section are Theorem \ref{theo: NS(Y) iff NS(X)}, where we determine the relation between $\NS(X)$ and $\NS(Y)$, and Theorem \ref{theo: theorem splitting in eigenspaces}. In the latter we determine the big and nef divisors on $Y$ which are associated by the map $\pi^*$ to specific ample divisors on $X$. This gives the relation between the dimension of the projective space in which we are embedding $X$ and the ones were $X/\langle\sigma\rangle$ has natural models. Examples of the geometric application of this result are provided in the next section.\\

We first fix an embedding of $\NS(X)$ in $H^2(X,\Z)\simeq \Lambda_{K3}$ and then we apply the map $\pi_*$ described in Proposition \ref{prop: pi_*} in order to find $\NS(Y)$.
To embed $\NS(X)$ in $\Lambda_{K3}$, we consider a specific embedding of $K_{12}$ in $\Lambda_{K3}$ constructed by embedding $\widetilde{K_{12}}$ in $U\oplus A_2(-1)\oplus E_6\oplus E_6\oplus E_6$ and then extending this embedding to the overlattice $K_{12}$ and $\Lambda_{K3}$.

We fixed a basis $\{k_i\}_{i=1,\ldots, 12}$ of the lattice $\widetilde{K_{12}}$ on which the bilinear form is the one given by the matrix \eqref{eq: K12tilde}. Then
$$\begin{array}{cccccccccc}
	\widetilde{K_{12}}&\stackrel{\lambda}{\rightarrow}&U\oplus&A_2(-1)\oplus &E_6\oplus&E_6\oplus &E_6\\  
k_i&\mapsto&(\underline{0},& \underline{0},&e_i^{(1)},&-e_i^{(2)},&\underline{0}),& \mbox{if }i=1,\ldots ,6\\
k_i&\mapsto&(\underline{0},& \underline{0},&e_i^{(1)},&\underline{0},&-e_i^{(3)}),& \mbox{if }i=7,\ldots ,12\end{array}$$
is an embedding of lattices.  
It extends to an embedding of $K_{12}$ into $\Lambda_{K3}$, indeed the overlattice $K_{12}$ is obtained by adding to $\lambda(\widetilde{K_{12}})$ the class $\lambda(z)$, written in \eqref{eq: n extra classes K12}.
This is a class contained also in the overlattice $\Lambda_{K3}$ of $U\oplus A_2(-1)\oplus E_6^{\oplus 3}$, so the enlargements of lattices are compatible.

This gives an embedding, still denoted by $\lambda$, of $K_{12}$ in $\Lambda_{K3}$, which is the $\mathbb{Q}$-linear extension of the embedding $\lambda$ defined above. We recall that the embedding of $K_{12}$ into $\Lambda_{K3}$ is unique up to isometries.

Both the lattices $\langle 2d\rangle\oplus K_{12}$ and $(\langle 2d\rangle\oplus K_{12})'$ admit a unique primitive embedding in $\Lambda_{K3}$ up to isometries and in the following proposition we exhibit one possible choice.
\begin{proposition}\label{prop: embedding NS(X) in LambdaK3} 
	The embedding 
$$\begin{array}{cccccccccc}
	\langle 2d\rangle&\stackrel{j}{\rightarrow}&U\oplus&A_2(-1)\oplus &E_6\oplus&E_6\oplus &E_6\\  
	L&\mapsto&(\left(\begin{array}{c}1\\d\end{array}\right),& \underline{0},&\underline{0},&\underline{0},&\underline{0})\end{array}$$
is such that $(j,\lambda): \langle 2d\rangle\oplus K_{12}\ra \Lambda_{K3}$ is a primitive embedding of $\langle 2d\rangle\oplus K_{12}$ in $\Lambda_{K3}$.

If $d\equiv 0\mod 3$, we consider the embedding 
$$\begin{array}{cccccccccc}
	\langle 2d\rangle&\stackrel{\widetilde{j}}{\rightarrow}&U\oplus&A_2(-1)\oplus &E_6\oplus&E_6\oplus &E_6\\  
	L&\mapsto&(\left(\begin{array}{c}3\\3k\end{array}\right),& \underline{0},&\underline{g}^{(1)},&\underline{g}^{(2)},&\underline{g}^{(3)}),\end{array}$$
where
$$\underline{g}^{(i)}=\left\{\begin{array}{ll}e_1^{(i)}+e_3^{(i)}+e_5^{(i)}&\mbox{ if }d\equiv 0\mod 9\\
e_1^{(i)}+e_3^{(i)}&\mbox{ if }d\equiv 3\mod 9\\
e_1^{(i)}&\mbox{ if }d\equiv 6\mod 9\\
\end{array}\right.\mbox{ and }k\mbox{ is s.t. }d=\left\{\begin{array}{ll}9k-9&\mbox{ if }d\equiv 0\mod 9\\
9k-3&\mbox{ if }d\equiv 3\mod 9\\
9k-6&\mbox{ if }d\equiv 6\mod 9\\
\end{array}\right.$$
Then $(\widetilde{j},\lambda): \langle 2d\rangle\oplus K_{12}\ra \Lambda_{K3}$ is an embedding whose primitive closure is $\left(\langle 2d \rangle\oplus K_{12}\right)'$.
\end{proposition}
\proof It is straightforward that $(j,\lambda)$ and $(\widetilde{j},\lambda)$ are embeddings of $\langle 2d\rangle\oplus K_{12}$ in $\Lambda_{K3}$. One can directly check that the first one is primitive, for example by observing that it coincides with the orthogonal complement of its orthogonal complement. The embeddings $\widetilde{j}:\langle 2d\rangle\ra\Lambda_{K3}$ and  $\lambda:K_{12}\ra\Lambda_{K3}$ are primitive. Nevertheless, $\widetilde{j}(L)+ (g^{(1)}- g^{(2)})+(g^{(1)} -g^{(3)})\in (\widetilde{j},\lambda)(\langle 2d\rangle\oplus K_{12})$ and it is divisible by 3 in $\Lambda_{K3}$. So the primitive closure of $(\widetilde{j},\lambda)(\langle 2d\rangle\oplus K_{12})$ is an overlattice of index 3 of $(\widetilde{j},\lambda)(\langle 2d\rangle\oplus K_{12})$ which contains primitively $K_{12}$ and $\langle 2d\rangle$, and thus it is $\left(\langle 2d \rangle\oplus K_{12}\right)'$.\endproof

\begin{theorem}\label{theo: NS(Y) iff NS(X)}
	Let $X$ be a projective K3 surface admitting a symplectic automorphism $\sigma$ of order 3 and let $Y$ be the resolution of $X/\langle\sigma\rangle$. Let $\rho(X)=13$.
	Then:
	
	$\NS(X)\simeq \langle 2d\rangle \oplus K_{12}$ if and only if $\NS(Y)\simeq \left(\langle 6d\rangle\oplus M_{\Z/3\Z}\right)';$ 
	
	$\NS(X)\simeq \left(\langle 6e\rangle\oplus K_{12}\right)'$ if and only if $\NS(Y)\simeq\langle 2e\rangle\oplus M_{\Z/3\Z}$. 
\end{theorem}
\proof To show the theorem, we apply $\pi_*$ to $\NS(X)=(j,\lambda)(\langle 2d\rangle, K_{12})\subset \Lambda_{K3}$ (resp. $(\widetilde{j},\lambda)(\langle 2d\rangle, K_{12})$). By Proposition \ref{prop: pi_*}, $$\pi_*(\lambda(K_{12}))=\underline{0}\in H^2(Y,\mathbb{Z})\simeq \Lambda_{K3}.$$ 

It remains to consider the image of $L$ under the embeddings we are considering.
In particular we have: 
$$\pi_*(j(L))=\left(\begin{array}{c}1\\d\end{array}\right)\subset U(3)\subset  H^2(Y,\mathbb{Z}).$$ Denoted by $H=\pi_*(j(L))$, we have $H^2=6d$. When we extend $U(3)\oplus A_2(-3)\oplus E_6$ to $H^2(Y,\mathbb{Z})$, the class $u_1'+du_2'$ glues with elements in $M_{\mathbb{Z}/3\mathbb{Z}}$ and gives the element $n_3+dn_4\in H^2(Y,\mathbb{Z})$, see Proposition \ref{prop: H^2(Y) extra classes}. Moreover $3(n_3+dn_4)\in \NS(Y)$ because it is a linear combination of classes in $\NS(Y)$, that is of $H=\pi_*(j(L))$ and of classes in $M_{\mathbb{Z}/3\mathbb{Z}}$. Since $n_3+dn_4\in H^2(Y,\mathbb{Z})$, $3\left(n_3+dn_4\right)\in \NS(Y)$ and $\NS(Y)$ is primitively embedded in $H^2(Y,\mathbb{Z})$, it follows that $n_3+dn_4\in \NS(Y)$.
Hence $\NS(Y)$ is spanned, over $\mathbb{Z}$, by $H$, the generators of $M_{\Z/3\Z}$ and an extra class, i.e. $n_3+dn_4$. So $\NS(Y)$ is an overlattice of index 3 of $\langle 6d\rangle\oplus M_{\mathbb{Z}/3\mathbb{Z}}$ and in particular it is necessarily the unique overlattice of index 3 described in Proposition \ref{prop:overlattice 2d+M3}.

We consider now the embedding $\widetilde{j}$: 
$$\pi_*(\widetilde{j}(L))=\left(\left(\begin{array}{c}3\\3k\end{array}\right),\underline{0}, 3\underline{g}\right)\subset U(3)\oplus A_2(-1)\oplus E_6\subset H^2(Y,\mathbb{Z}),$$
where $\underline{g}$ is the vector $\underline{g}^{(i)}$ of Proposition \ref{prop: embedding NS(X) in LambdaK3}. It is clear that $\pi_*(\widetilde{j}(L))$ is not primitive, since it is 3 divisible. So $$\pi_*(\widetilde{j}(L))/3=\left(\left(\begin{array}{c}1\\k\end{array}\right),\underline{0}, \underline{g}\right)\in \NS(Y)\subset H^2(Y,\mathbb{Z}).$$
and we define $H$ to be $\pi_*(L)/3$. So $H^2=6k-(\underline{g})^2$. Even enlarging $U(3)\oplus A_2(-3)\oplus E_6$ to $H^2(Y)$ by gluing the classes of $M_{\mathbb{Z}/3\mathbb{Z}}$ we do not find new classes in $\NS(Y)$, so $\NS(Y)$ is generated, over $\mathbb{Z}$, by $H$ and by the classes generating $M_{\mathbb{Z}/3\mathbb{Z}}$. The intersection form is $\langle 2d/3\rangle\oplus M_{\mathbb{Z}/3\mathbb{Z}}$, which concludes the proof.\endproof

\begin{corollary}
The K3 surface $X$ is a generic member of the family of the $\left(\langle 2d\rangle\oplus K_{12}\right)$-polarized K3 surfaces if and only if $Y$	 is a generic member of the family of the $\left(\langle 2d\rangle\oplus M_{\Z/3\Z}\right)'$-polarized K3 surfaces.

The K3 surface $X$ is a generic member of the family of the $\left(\langle 6e\rangle\oplus K_{12}\right)'$-polarized K3 surfaces if and only if $Y$ is a generic member of the family of the $\left(\langle 2e\rangle\oplus M_{\Z/3\Z}\right)$-polarized K3 surfaces.
\end{corollary}


To apply the map $\pi^*$ to $\NS(Y)$, we first fix an embedding of $M_{\Z/3\Z}$ in $\Lambda_{K3}$:
we constructed $\Lambda_{K3}\simeq H^2(Y,\Z)$ as overlattice of index $3^4$ of $A_2(-1)\oplus U(3)\oplus E_6\oplus M_{\Z/3\Z}$ in Proposition \ref{prop: H^2(Y) extra classes}. The natural embedding of $M_{\Z/3\Z}$ in $A_2(-1)\oplus U\oplus E_6\oplus M_{\Z/3\Z}$,
$$\begin{array}{cccccccccc}\mu:&M_{\Z/3\Z}&\ra&A_2(-1)\oplus& U(3)\oplus &E_6\oplus &M_{\Z/3\Z}&\hookrightarrow \Lambda_{K3}\\&m_i&\mapsto&(\underline{0},&\underline{0},&\underline{0},&m_i)\end{array}$$ extends to a primitive embedding $\mu:M_{\Z/3\Z}\ra \Lambda_{K3}$.

Both the lattices $\langle 2e\rangle\oplus M_{\Z/3\Z}$ and $(\langle 2e\rangle\oplus M_{\Z/3\Z})'$ admit a unique primitive embedding in $\Lambda_{K3}$ up to isometries and in the following proposition we exhibit one possible choice.
\begin{proposition}\label{prop: embedding of NS(Y) in LambdaK3} 
	The embedding 
	$$\begin{array}{cccccccccc}
	\langle 2e\rangle&\stackrel{h}{\rightarrow}&U(3)\oplus&A_2(-1)\oplus &E_6\oplus&M_{\Z/3\Z}\\  
	H&\mapsto&(\left(\begin{array}{c}1\\k\end{array}\right),& \underline{0},&\underline{f},&\underline{0})\end{array}$$
	where $$\underline{f}=\left\{\begin{array}{lll
}e_1+e_3+e_5&\mbox{ if }&2e=6k-6,\\
	e_1+e_3&\mbox{ if }&2e=6k-4,\\
	e_1&\mbox{ if }&2e=6k-2\\
	\end{array}\right.$$
	
	is such that $(h,\mu): \langle 2e\rangle\oplus M_{\Z/3\Z}\ra \Lambda_{K3}$ is a primitive embedding of $\langle 2e\rangle\oplus M_{\Z/3\Z}$ in $\Lambda_{K3}$.
	
	If $e\equiv 0\mod 3$, then there exists $k$ such that $e=3k$ and the embedding 
	$$\begin{array}{cccccccccc}
	\langle 2e\rangle&\stackrel{\widetilde{h}}{\rightarrow}&U(3)\oplus&A_2(-1)\oplus &E_6\oplus&M_{\Z/3\Z}\\  
	H&\mapsto&(\left(\begin{array}{c}1\\k\end{array}\right),& \underline{0},&\underline{0},&\underline{0}),\end{array}$$
	is such that $(\widetilde{h},\lambda): \langle 2e\rangle\oplus M_{\Z/3\Z}\ra \Lambda_{K3}$ is an embedding whose primitive closure is $\left(\langle 2e \rangle\oplus M_{\Z/3\Z}\right)'$.
\end{proposition}
\proof The proof is analogous to the one of Proposition \ref{prop: embedding NS(X) in LambdaK3}.\endproof


\begin{corollary}\label{cor: pullback Di}
	Let $Y$ be a K3 surface such that $\NS(Y)\simeq \langle 2e\rangle\oplus M_{\Z/3\Z}$ and the embedding of $\NS(Y)$ in $\Lambda_{K3}$ is $(h,\mu)$.
	Then $$\pi^*((h,\mu)(D_i))=\widetilde{j}(L)$$ where $D_i$, $i=1,2,3$ are the divisors defined in \eqref{eq: D_i not tilde}.
	
	Let $Y$ be a K3 surface such that $\NS(Y)\simeq (\langle 2e\rangle\oplus M_{\Z/3\Z})'$ and the embedding of $\langle 2e\rangle\oplus M_{\Z/3\Z}$ in $\Lambda_{K3}$ is $(\widetilde{h},\mu)$.
	Then $$\pi^*((\widetilde{h},\mu)(D_i))=j(L)$$ where $D_i$, $i=1,2,3$ are the divisors defined in \eqref{eq: D_i case tilde 1}, \eqref{eq: D_i case tilde 2}, \eqref{eq: D_i case tilde 3}.
\end{corollary} 
\proof
We first consider $\pi^*$ of $h(H)$ and of $\widetilde{h}(H)$:
$$\pi^*(h(H))=\pi^*\left(\left(\begin{array}{c}1\\k\end{array}\right), \underline{0},\underline{f},\underline{0}\right)=\left(\left(\begin{array}{c}3\\3k\end{array}\right), \underline{0},\underline{g}^{(1)},\underline{g}^{(2)},\underline{g}^{(3)}\right)\in U\oplus A_2(-1)\oplus E_6\oplus E_6\oplus E_6$$

and 
$$\pi^*(\widetilde{h}(H))=\pi^*\left(\left(\begin{array}{c}1\\k\end{array}\right), \underline{0},\underline{0},\underline{0}\right)=\left(\left(\begin{array}{c}3\\3k\end{array}\right), \underline{0},\underline{0},\underline{0},\underline{0},\right)\in U\oplus A_2(-1)\oplus E_6\oplus E_6\oplus E_6.$$
In particular we observe that $$\pi^*(h(H))=\widetilde{j}(L),\ \ \pi^*(\widetilde{h}(H))=3j(L).$$

The divisors $D_i$ defined in Remark \ref{rem: defi Di} are embedded in $\NS(Y)\subset \Lambda_{K3}$ by the embedding $(h,\mu)$ or $(\widetilde{h},\mu)$ (according to the properties of $\NS(Y)$). We observe that $\pi^*(M_i)=\underline{0}\in H^2(X,\Z)$, which allows to conclude. As example we consider the divisor $D_1$ described in Remark \ref{rem: defi Di}, \eqref{eq: D_i case tilde 1}: $$\pi^*((\widetilde{h},\mu)(D_1))=\pi^*\left(\frac{\widetilde{h}(H)-\mu(M_1^{(1)}+2M_2^{(1)}+M_1^{(2)}+2M_2^{(2)}+M_1^{(3)}+2M_2^{(3)})}{3}\right)=j(L).$$ \endproof

\begin{theorem}\label{theo: theorem splitting in eigenspaces}
	Let $X$ be a projective  K3 surface admitting a symplectic automorphism of order 3 such that $\rho(X)=13$, $L$ be the ample generator of $K_{12}^{\perp_{NS(X)}}$, $Y$ be the minimal resolution of $X/\langle\sigma\rangle$, and $D_i$, $i=1,2,3$ the divisors defined in Remark \eqref{rem: defi Di}. Then 
	$$H^0(X,L)=\pi^*H^0(Y,D_1)\oplus \pi^*H^0(Y,D_2)\oplus \pi^*H^0(Y,D_3)$$
	and the previous decomposition corresponds to the decomposition of $H^0(X, L)$ in eigenspaces with respect to the action of $\sigma^*$ on $H^0(X, L)$.
	\end{theorem}
\proof We recall that the action of $\sigma^*$ on the divisor $L$ is the identity (since  $K_{12}\simeq (\NS(X)^{\sigma^*})^{\perp}$). Hence $\sigma^*$ acts on the vector space $V:=H^0(X,L)$ and its action splits $V$ in the direct sum of three eigenspaces, i.e. $V:=V_{+1}\oplus V_{\zeta_3}\oplus V_{\zeta_3^2}$.  

By Corollary \ref{cor: pullback Di}, the pullbacks of the sections in $H^0(Y,D_i)$ are sections in $H^0(X,L)$. Moreover, $D_i$ are divisors on $Y$, so their sections are well defined on the quotient surface $X/\langle\sigma\rangle$. Hence, given a basis $\{s_1,\ldots s_r\}$ of $H^0(Y,D_i)$, $\pi^*(s_1),\ldots \pi^*(s_r)$ lie in the same eigenspace for the action of $\sigma^*$ on $H^0(X,L)$, otherwise they would not be contained in the same $H^0(Y,B)$ for a divisor $B\in \NS(Y)$. The images of the exceptional curves $M_i^{(k)}$ under the maps $\varphi_{|D_i|}:Y\ra\mathbb{P}(H^0(Y,D_i)^{\vee})$ changes if one considers the divisors $D_i$ or $D_j$ with $i\neq j$ (since the intersection number of the divisors $M_h^{(k)}$ with $D_i$ and $D_j$ are not the same). So the sections of $D_i$ and $D_j$ lie on different eigenspaces $V_{\epsilon}\subset V$ where $\epsilon=+1,\zeta_3,\zeta_3^2$. It remains to prove that $\pi^*H^0(Y,D_1)$ coincides with one eigenspace $V_{\epsilon}$ and it is not only contained in it. To this purpose, it suffices to prove that $$\dim(H^0(X,L))=\dim(H^0(Y,D_1))+\dim(H^0(Y,D_2))+\dim(H^0(Y,D_3)).$$ We already proved that $\pi^*H^0(Y,D_1)$ is contained in an eigenspace, so we have $$\dim(H^0(X,L))\geq \dim(H^0(Y,D_1))+\dim(H^0(Y,D_2))+\dim(H^0(Y,D_3))$$ and it remains to prove that $$\dim(H^0(X,L))\leq \dim(H^0(Y,D_1))+\dim(H^0(Y,D_2))+\dim(H^0(Y,D_3)).$$

If $\NS(X)\simeq \langle 2d\rangle\oplus K_{12}$ and $L$ is the ample generator of $\langle 2d\rangle$, then $\dim(H^0(X,L))=d+2$, $\NS(Y)\simeq \left(\langle 6d\rangle\oplus M_{\Z/3\Z}\right)'$, and by Riemann--Roch theorem we obtain:
\begin{itemize}
	\item if $3d\equiv 0\mod 9$ then $\chi(D_1)=\frac{d}{3}+1$, $\chi(D_2)=\frac{d}{3}+1$, $\chi(D_3)=\frac{d}{3}$;
	\item if $3d\equiv 3\mod 9$ then $\chi(D_1)=\frac{d}{3}+\frac{2}{3}$, $\chi(D_2)=\frac{d}{3}+\frac{2}{3}$, $\chi(D_3)=\frac{d}{3}+\frac{2}{3}$;  
	\item if $3d\equiv 6\mod 9$ then $\chi(D_1)=\frac{d}{3}+\frac{4}{3}$, $\chi(D_2)=\frac{d}{3}+\frac{1}{3}$, $\chi(D_3)=\frac{d}{3}+\frac{1}{3}$.
\end{itemize}
In all the listed cases $\chi(D_1)+\chi(D_2)+\chi(D_3)=d+2=\dim(H^0(X,L))$.

The divisors $D_i$ have a positive intersection with the pseudoample divisor $H$ and their self intersection is bigger than or equal to $-2$. Then $h^0(D_i)>0$ and $h^2(D_i)=h^0(-D_i)=0$. Hence $h^0(D_i)\geq \chi (D_i)$. 
So $$h^0(Y,D_1)+h^0(Y,D_2)+h^0(Y,D_3)\geq \chi(D_1)+\chi(D_2)+\chi(D_3)=h^0(X,L)=d+2,$$ which concludes the proof if $\NS(Y)=\langle 2d\rangle\oplus K_{12}$.

If $\NS(X)\simeq \left(\langle 2d\rangle\oplus K_{12}\right)'$ and $L$ is the ample generator of $\langle 2d\rangle$, then we argue as above, observing that $\dim(H^0(X,L))=d+2$, $\NS(Y)\simeq \langle \frac{2d}{3}\rangle\oplus M_{\Z/3\Z}$ and
$\chi(D_1)=\frac{d}{3}+2$, $\chi(D_2)=\frac{d}{3}$, $\chi(D_3)=\frac{d}{3}$. \endproof

\section{Examples}\label{sec: examples}

We saw in Theorem \ref{theor:NS X} that a K3 surface $X$ such that $\NS(X)$ is either $\langle 2d\rangle\oplus K_{12}$ or $\left(\langle 2d\rangle\oplus K_{12}\right)'$ is projective and admits an order 3 symplectic automorphism $\sigma$. In Theorem \ref{theo: theorem splitting in eigenspaces} we observed that $\sigma^*$ acts on the space $H^0(X, L)^{\vee}$, where $L$ is the ample generator of $\langle 2d\rangle$ in the N\'eron--Severi group. Hence $X$ admits a projective model in  $\mathbb{P}(H^0(X, L)^{\vee})$ such that the automorphism $\sigma$ is the restriction to $\varphi_{|L|}(X)$ of a projective transformation of $\mathbb{P}(H^0(X, L)^{\vee})$. As an application of the Theorem \ref{theo: theorem splitting in eigenspaces} we can determine this projective transformation on $\mathbb{P}(H^0(X, L)^{\vee})$, in particular the dimension of its eigenspaces.

The surface $Y$, minimal resolution of $X/\langle\sigma\rangle$,
has a polarization $H$ induced by $L$ and described in Theorem \ref{theo: NS(Y) iff NS(X)}. It is orthogonal to all the classes $M_i^{(j)}$, so $\varphi_{|H|}(Y)$ is the singular model of $Y$ where all the curves $M_{i}^{(j)}$ are contracted, thus it is the projective model of $X/\langle\sigma\rangle$.

As an application of the Theorem \ref{theo: NS(Y) iff NS(X)} we have a relation between the dimension of the projective space $\mathbb{P}(H^0(X,L)^{\vee})$, where the surface $X$ is embedded, and the dimension of the projective space $\mathbb{P}(H^0(Y,H)^{\vee})$, where the quotient surface $X/\langle\sigma\rangle$ is embedded. The aim of this section is to apply the previous results to obtain explicit projective equations of $X$ and $X/\langle\sigma\rangle$ for certain values of $d$.

Moreover, by Theorem \ref{theo: theorem splitting in eigenspaces} we find other projective models of the surface $Y$, related with the existence of the divisors $D_i$. 

\subsection{The case $d=1$} In this case there is a unique possibility for $\NS(X)$, which is $\NS(X)=\langle 2\rangle\oplus K_{12}$; $\varphi_{|L|}:X\ra\mathbb{P}^2=\mathbb{P}(H^0(X,L)^{\vee})$ and $\sigma$ is induced by an automorphism of $\mathbb{P}^2$. 
By Theorem \ref{theo: NS(Y) iff NS(X)}, $\NS(Y)\simeq (\langle 6\rangle\oplus M_{\Z/3\Z})'$ and the definition of the divisors $D_i$ on $Y$ is the one given in Remark \ref{rem: defi Di}, \eqref{eq: D_i case tilde 2}. So $\dim(H^0(Y,D_i))=1$, for $i=1,2,3$, and by Theorem \ref{theo: theorem splitting in eigenspaces}, the eigenspaces $V_{\epsilon}$ have dimension 1.

Therefore there exists a choice of coordinates of $\mathbb{P}^2=\mathbb{P}(H^0(X,L)^{\vee})$ such that the action of $\sigma$ on $\mathbb{P}^2$ is 
$$\sigma:\mathbb{P}^2\ra\mathbb{P}^2\ \ \ (x_0:x_1:x_2)\mapsto(x_0:\zeta_3 x_1:\zeta_3^2 x_2).$$
Let us consider the equation $f_6$ of the plane sextic curves invariant for $\sigma$: $$a_1x_0^6+a_2x_0^4x_1x_2+a_3x_0^3x_1^3+a_4x_0^3x_2^3+a_5x_0^2x_1^2x_2^2+a_6x_0x_1^4x_2+a_7x_0x_1x_2^4+a_8x_1^6+a_9x_1^3x_2^3+a_{10}x_2^6.$$
The double cover of $\mathbb{P}^2$ branched on these sextic curves has equation $$w^2=f_6(x_0:x_1:x_2)\subset \mathbb{P}(3,1,1,1).$$ This gives a family of K3 surfaces $X$ admitting a symplectic automorphism of order 3, which is a lift of $\sigma$. 
The family depends on 10 parameters, but it is defined up to projective transformations of $\mathbb{P}^2$ which commute with $\sigma$. These are the diagonal transformations, so the dimension of the family, up to projective transformations is $(10-1)-(3-1)=7$, which is the expected dimension of a family of projective K3 surfaces with N\'eron--Severi group isometric to $\langle 2\rangle\oplus K_{12}$. The points of $\mathbb{P}^2$ fixed by $\sigma$ are $(1:0:0)$, $(0:1:0)$, $(0:0:1)$. None of them lies on the branch sextic curve, hence they correspond to 6 points on the K3 surface $X$ and there exists a lift of $\sigma$ which fixes all these points. So the automorphism $\sigma$ acts on the K3 surface fixing 6 isolated points, and hence it is symplectic, therefore we found an explicit equation of the family of K3 surfaces whose N\'eron--Severi group is $\langle 2\rangle\oplus K_{12}$.

By Theorem \ref{theo: NS(Y) iff NS(X)}, the divisor $H$ on $Y$ has self intersection 6 and then $\varphi_{|H|}(Y)\subset\mathbb{P}^4$. The following functions are invariant for $\sigma$
$$z_0:=w, z_1:=x_0^3, z_2:=x_1^3, z_3:=x_2^3, z_4:=x_0x_1x_2$$
and satisfies the equations 
\begin{equation}\label{eq: equation of Y case d=1}\left\{\begin{array}{l}
z_1z_2z_3=z_4^3\\
z_0^2=a_1z_1^2+a_2z_1z_4+a_3z_1z_2+a_4z_1z_3+a_5z_4^2+a_6z_2z_4+a_7z_3z_4+a_8z_2^2+a_9z_2z_3+a_{10}z_3^2.\end{array}\right.\end{equation}
Hence \eqref{eq: equation of Y case d=1} are the equations of $X/\langle\sigma\rangle$ in $\mathbb{P}^4$, i.e. the equations of $\varphi_{|H|}(Y)$.

\subsection{The case $d=2$} In this case there is a unique possibility for $\NS(X)$, which is $\NS(X)=\langle 4\rangle\oplus K_{12}$,  $\varphi_{|L|}:X\ra\mathbb{P}^3=\mathbb{P}(H^0(X,L)^{\vee})$ and $\sigma$ is induced by an automorphism of $\mathbb{P}^3$. 
By Theorem \ref{theo: NS(Y) iff NS(X)}, $\NS(Y)\simeq (\langle 12\rangle\oplus M_{\Z/3\Z})'$ and the definition of the divisors $D_i$ on $Y$ is the one given in Remark \ref{rem: defi Di}, \eqref{eq: D_i case tilde 3}. So $\dim(H^0(Y,D_1))=2$, $\dim(H^0(Y, D_i))=1$ for $i=2,3$ (cf. proof of Theorem \ref{theo: theorem splitting in eigenspaces}) hence one eigenspace has dimension 2 and the other two have dimension 1.

So there exists a choice of coordinates of $\mathbb{P}^3=\mathbb{P}(H^0(X,L)^{\vee})$ such that the action of $\sigma$ on $\mathbb{P}^3$ is 
$$\sigma:\mathbb{P}^3\ra\mathbb{P}^3\ \ \ (x_0:x_1:x_2:x_3)\mapsto(x_0:x_1:\zeta_3 x_2:\zeta_3^2 x_3).$$
The quartic equations invariant for $\sigma$ are 
 \begin{align}\label{eqn:quartic}f_4(x_0:x_1)+f_2(x_0:x_1)x_2x_3+f_1(x_0:x_1)x_2^3+g_1(x_0:x_1)x_3^3+\alpha x_2^2x_3^2 =0\end{align}
 where $f_i$ and $g_i$ are homogeneous polynomials of degree $i$.
 
This defines a family of K3 surfaces admitting an automorphism induced by $\sigma$.
The family depends on 13 parameters but it is defined up to the actions of the projective transformation of $\mathbb{P}^3$ which commutes with $\sigma$. So the family depends on $(13-1)-(6-1)=7$ parameters.
There are 4 fixed points in the eigenspace $V_{+1}$, which are defined by $f_4(x_0:x_1)=0,\ x_2=x_3=0$; 1 fixed point in $V_{\zeta_3}$, i.e. $(0:0:1:0)$, and 1 in $V_{\zeta_3^2}$, i.e. $(0:0:0:1)$. In particular the automorphism $\sigma$ fixes 6 isolated points on the K3 surfaces in the family and hence it is a symplectic automorphism of each member of this family.


By Theorem \ref{theo: NS(Y) iff NS(X)}, the divisor $H$ has self intersection 12 and then $\varphi_{|H|}(Y)\subset\mathbb{P}^{7}$. We now show that the ideal defining $\varphi_{|H|}(Y)\subset\mathbb{P}^{7}$ is generated by 10 quadrics and we determine a set of generators.

The N\'eron--Severi group $\NS(Y)$ is isometric to $(\langle 12\rangle\oplus M_{\Z/3\Z})'$ and hence each class in $\NS(Y)$ can be written as $\alpha H+\beta m$ where $m\in M_{\Z/3\Z}$ and $\alpha,\ \beta\in \frac{1}{3}\Z$. Let $C$ be an irreducible curve on $Y$, hence $HC\geq 0 $ and so the class of $C\in\NS(Y)$ is $\alpha H+\beta m$ with $\alpha\geq 0$. The intersection $HC$ is either 0, if $\alpha=0$ or $12\alpha\geq 4$. In particular there are no curves $C$ such that $HC=2$. This implies, by \cite[Theorem 5.2]{SD}, that the linear system $H$ is not hyperelliptic and by \cite[Theorem 7.2]{SD}, that the ideal of $\varphi_{|H|}(Y)\subset \mathbb{P}^7$ is generated by quadrics. Since $h^0(Y,H)=10$, $\dim S^2(H^0(Y,H))=\binom{9}{2}=36$ and, by $(2H)^2=28$, it follows $h^0(Y,2H)=26$. Hence the ideal of $\varphi_{|H|}(Y)\subset \mathbb{P}^7$ is generated by 10 quadrics and we are going to determine them.

We preliminary observe that the number $n_h$ of monomials of degree $h$ in the variables $x_2$ and $x_3$ which are invariant for the action of $\sigma$ are as in the following table
$$\begin{array}{|c||c|c|c|c|c|c|c|}

\hline
h&0&1&2&3&4&5&6\\
\hline
n_h&1&0&1&2&1&2&3\\
\hline\end{array}$$
This allows to compute the number $m_k$ of monomials of degree $k$ in the variables $x_i$, $i=0,1,2,3$ which are invariant for the action of $\sigma$ are as in the following table
$$\begin{array}{|c||c|c|c|c|c|c|c|}
		\hline
	k&0&1&2&3&4&5&6\\
	\hline
	m_k&1&2&4&8&13&22&30\\
	\hline\end{array}.$$

There are 8 invariant monomials of degree 3 in the variables $x_i$:
$$a_i=x_i^3, i=0,1,2,3,\ \ a_4=x_0^2x_1,\ a_5=x_0x_1^2,\ a_6=x_0x_2x_3,\ a_7=x_1x_2x_4.$$
We choose them as coordinates of the projective space $\mathbb{P}^7_{a_i}$
such that $\varphi_{|H|}(Y)\subset\mathbb{P}^7_{a_i}$.

The following quadric equations are satisfied by the variables $a_i$:
\begin{align}\label{6 quadrics}a_4^2=a_0a_5,\ a_4a_5=a_0a_1,\ a_4a_6=a_0a_7,\ a_5^2=a_1a_4, \ a_5a_6=a_4a_7,\ a_5a_7=a_1a_6\end{align}
so they are contained in the ideal of quadrics defining $\varphi_{|H|}(Y)\subset\mathbb{P}^7_{a_i}$.

We now determine the other 4 quadrics generating the ideals of $\varphi_{|H|}(Y)$.

The 4 monomials in $x_i$ which are of degree 2 and are invariant for $\sigma$, are $x_0^2$, $x_1^2$, $x_0x_1$, $x_2x_3$.
Since the quartic \eqref{eqn:quartic} is invariant for the action of $\sigma$, multiplying it for each of the invariant monomials of degree 2 written above one finds an invariant sextic in the variables $x_i$, i.e. the followings, with $i=1,2$ and $(h,k)=(0,1), (2,3)$:
\begin{align}\label{eqn: sextic}\begin{array}{c}x_i^2\left(f_4(x_0:x_1)+f_2(x_0:x_1)x_2x_3+f_1(x_0:x_1)x_2^3+g_1(x_0:x_1)x_3^3+\alpha x_2^2x_3^2\right)\\
x_hx_k\left(f_4(x_0:x_1)+f_2(x_0:x_1)x_2x_3+f_1(x_0:x_1)x_2^3+g_1(x_0:x_1)x_3^3+\alpha x_2^2x_3^2\right).\end{array}\end{align}

The monomials appearing in these sextics are invariant monomials of degree 6. Each of them can be expressed as a monomial of degree 2 in the $a_i$'s, since the space of the quadric in the variables $a_i$ has dimension 30 (there 36 quadrics and 6 relations, as seen above) which is the same dimension of the space of the invariant polynomials of degree 6 in the $x_i$'s.

So each of the four sextics in \eqref{eqn: sextic} can be expressed as a quadric equation in the variables $a_i$ and this provides other four quadrics in $\mathbb{P}^7$ which vanish one the surface $\varphi_{|H|}(Y)$.

These four quadrics together with the ones listed in \eqref{6 quadrics} are the 10 generators of the ideal of $\varphi_{|H|}(Y)$.

The map induced by $D_1$ is an elliptic fibration on the quotient surface, i.e. $\varphi_{|D_1|}:Y\ra\mathbb{P}^1$ is an elliptic fibration given by the projection of the quartic surface from the space $x_0=x_1=0$ to the linear subspace $\mathbb{P}^1_{(x_0:x_1)}\subset\mathbb{P}^3_{(x_0:x_1:x_2:x_3)}$.

\subsection{The case $d=3$}\label{subsec: d=3} In this case there are two possibilities for $\NS(X)$, which are $\NS(X)=\langle 6\rangle\oplus K_{12}$ and $\NS(X)=(\langle 6\rangle\oplus K_{12})'$, in both cases $\varphi_{|L|}:X\ra\mathbb{P}^4=\mathbb{P}(H^0(X,L)^{\vee})$ and $\sigma$ is induced by an automorphism of $\mathbb{P}^4$.
\subsubsection{The case $\NS(X)=\langle 6\rangle\oplus K_{12}$.}  
By Theorem \ref{theo: NS(Y) iff NS(X)}, $\NS(Y)\simeq (\langle 18\rangle\oplus M_{\Z/3\Z})'$ and the definition of the divisors $D_i$ on $Y$ is the one given in Remark \ref{rem: defi Di}, \eqref{eq: D_i case tilde 1}. So $\dim(H^0(Y,D_i))=2$, $i=1,2$, $\dim(H^0(Y, D_3))=1$ (cf. proof of Theorem \ref{theo: theorem splitting in eigenspaces}) hence two eigenspaces have dimension 2 and the other one has dimension 1.

So there exists a choice of coordinates of $\mathbb{P}^4=\mathbb{P}(H^0(X,L)^{\vee})$ such that the action of $\sigma$ on $\mathbb{P}^4$ is 
$$\sigma:\mathbb{P}^4\ra\mathbb{P}^4\ \ \ (x_0:x_1:x_2:x_3:x_4)\mapsto(x_0:x_1:\zeta_3 x_2:\zeta_3 x_3:\zeta_3^2 x_4).$$
The K3 surfaces embedded in $\mathbb{P}^4$ are complete intersections of a quadric and a cubic, hence one has to find equations for cubic and quadric hypersurfaces such that both the hypersurfaces are invariant for the action of $\sigma$ and the restriction of $\sigma$ to their intersection fixes exactly 6 points. A choice which satisfies all the required conditions is:
$$\left\{\begin{array}{l}a_1(x_0:x_1)b_1(x_2:x_3)+\alpha x_4^2=0,\\
c_3(x_0:x_1)+d_1(x_0:x_1)e_1(x_2:x_3)x_4+f_3(x_2:x_3)+\beta x_4^3=0,\end{array}
\right.$$
where all the polynomials are homogeneous and their degrees are their subscripts.  
There are 3 fixed points in $x_0=x_1=x_4=0$ (the ones which satisfy $f_3(x_2:x_3)=0$) and other 3 fixed points in  $x_2=x_3=x_4=0$ (the ones which satisfy $c_3(x_0:x_1)=0$).

By Theorem \ref{theo: NS(Y) iff NS(X)}, the divisor $H$ has self intersection 18, and then $\varphi_{|H|}(Y)\subset\mathbb{P}^{10}$. The maps $\varphi_{|D_1|}:Y\ra\mathbb{P}^1$ and $\varphi_{|D_2|}:Y\ra\mathbb{P}^1$ define two elliptic fibrations, each of them has 3 independent sections corresponding to certain curves $M_i^{(j)}$ and contracts all the other curves $M_i^{(j)}$, which are necessarily irreducible components of reducible fibers. The elliptic fibrations correspond to the projections of the surface to the subspaces $\mathbb{P}^1_{(x_0:x_1)}$ and $\mathbb{P}^1_{(x_3:x_4)}$ respectively.

\subsubsection{The case $\NS(X)=\left(\langle 6\rangle\oplus K_{12}\right)'$.}  
By Theorem \ref{theo: NS(Y) iff NS(X)}, $\NS(Y)\simeq \langle 2\rangle\oplus M_{\Z/3\Z}$ and the definition of the divisors $D_i$ on $Y$ is the one given in Remark \ref{rem: defi Di}, equation \eqref{eq: D_i not tilde}. So $\dim(H^0(Y,D_1))=3$,  $\dim(H^0(Y, D_i))=1$, $i=2,3$ (cf. proof of Theorem \ref{theo: theorem splitting in eigenspaces}) hence one eigenspace has dimension 3 and the others have dimension 1.

So there exists a choice of coordinates of $\mathbb{P}^4=\mathbb{P}(H^0(X,L)^{\vee})$ such that the action of $\sigma$ is 
$$\sigma:\mathbb{P}^4\ra\mathbb{P}^4\ \ \ (x_0:x_1:x_2:x_3:x_4)\mapsto(x_0:x_1: x_2:\zeta_3 x_3:\zeta_3^2 x_4).$$
The K3 surface $X$ is the complete intersection of a quadric and a cubic whose equations are invariant for $\sigma$:
$$\left\{\begin{array}{l}q_2(x_0:x_1:x_2)+\alpha x_3x_4=0,\\
c_3(x_0:x_1:x_2)+l_1(x_0:x_1:x_2)x_3x_4+\beta x_3^3+\gamma x_4^3=0,\end{array}\right.$$
where all the polynomials are homogeneous and their degrees are the subscript numbers.  
The 6 fixed points are all contained in the eigenspace $V_{+1}$, which are the intersections of $c_3=0$ and $q_2=0$ in the plane $x_3=x_4=0$.

By Theorem \ref{theo: NS(Y) iff NS(X)}, the divisor $H$ has self intersection 2 and then $\varphi_{|H|}:Y\ra \mathbb{P}^{2}$ is a double cover. This map is induced by the projection of $X$ on the invariant plane $\mathbb{P}^2_{(x_0:x_1:x_2)}$. The projection to $\mathbb{P}^2_{(x_0:x_1:x_2)}$ is the projection from the line $x_0=x_1=x_2=0\subset \mathbb{P}^4$. Let us denote by $p:\mathbb{P}^4\ra\mathbb{P}^2$ this projection and by $p_X$ its restriction to $X$. 
Given a point $(\overline{x_0}:\overline{x_1}:\overline{x_2})\in\mathbb{P}^2$, its inverse image $p^{-1}_X((\overline{x_0}:\overline{x_1}:\overline{x_2}))$ consists of six points. Indeed in the affine chart $x_4=1$, the inverse image of the point for the map $p$ is the plane whose parametric equations are
$$\left\{\begin{array}{l}x_0=s \overline{x_0},\\x_1=s \overline{x_1},\\x_2=s \overline{x_2},\\x_3=t.\end{array}\right.$$
The intersection of this plane with $X$ is
$$\left\{\begin{array}{l}t=-s^2q_2(\overline{x_0}:\overline{x_1}:\overline{x_2})/\alpha,\\
s^6\beta \left(q_2(\overline{x_0}:\overline{x_1}:\overline{x_2})/\alpha\right)^3+s^3\left(c_3(\overline{x_0}:\overline{x_1}:\overline{x_2})-l_1(\overline{x_0}:\overline{x_1}:\overline{x_2})q_2(\overline{x_0}:\overline{x_1}:\overline{x_2})/\alpha\right)+\gamma=0.
\end{array}\right.$$
If $$\left(c_3(\overline{x_0}:\overline{x_1}:\overline{x_2})-\frac{l_1(\overline{x_0}:\overline{x_1}:\overline{x_2})q_2(\overline{x_0}:\overline{x_1}:\overline{x_2})}{\alpha}\right)^2+4\gamma\beta \frac{q_2(\overline{x_0}:\overline{x_1}:\overline{x_2})^3}{\alpha^3}\neq 0,$$ then one obtains two 2 solutions for $s^3$ in the last equation. So one obtains 6 solutions for $s$ and each choice for $s$ determines a choice for $t$ and thus a point in $X$.
Hence $p_x^{-1}(\overline{x_0}:\overline{x_1}:\overline{x_2})$ consists of 6 points. These six points form two orbits for the action of $\sigma$ (indeed $\sigma$ acts multiplying $s$ by a root of unity, so it identifies values of $s$ which have the same third power). Hence these six points on $X$ correspond to 2 points on the quotient surface $X/\langle\sigma\rangle$.

On the other hand if $$\left(c_3(\overline{x_0}:\overline{x_1}:\overline{x_2})-\frac{l_1(\overline{x_0}:\overline{x_1}:\overline{x_2})q_2(\overline{x_0}:\overline{x_1}:\overline{x_2})}{\alpha}\right)^2+4\beta \frac{q_2(\overline{x_0}:\overline{x_1}:\overline{x_2})^3}{\alpha^3}=0,$$ then one obtains one solution (with multiplicity 2) for $s^3$. Hence in this case $p_X^{-1}(\overline{x_0}:\overline{x_1}:\overline{x_2})$ consists of 3 points which are in the same orbit for $\sigma$ and determines a point in $X/\langle\sigma\rangle$.

So this projection describes $X/\langle\sigma\rangle$ as a double cover of $\mathbb{P}^2_{(x_0:x_1:x_2)}$ branched on the sextic  $$\alpha^3c_3^2+\alpha l_1^2q_2^2-2\alpha^2c_3l_1q_2+4\beta q_2^3=0$$
where we wrote $c_3$ (resp. $l_1$ and $q_2$) instead of $c_3(x_0:x_1:x_2)$ (resp. $l_1(x_0:x_1:x_2)$ and $q_2(x_0:x_1:x_2)$).
We observe that the branch locus is a sextic with 6 singularities of type $A_2$, that are the points where $c_3(x_0:x_1:x_2)=q_2(x_0:x_1:x_2)=0$. These are the images of the 6 points in $X$ fixed by $\sigma$.

\section{Applications to isogenies and generalized Shioda--Inose structures}\label{sec: applications}

In the previous sections, we described the action induced in cohomology by a symplectic automorphism of order 3, (see Section \ref{sec: cohomological action}) and the relation between the N\'eron--Severi group of projective K3 surfaces admitting an order 3 symplectic automorphism and the N\'eron--Severi group of its quotient (Theorem \ref{theo: NS(Y) iff NS(X)}). Hence, we are able to generalize some of the results obtained for the involutions to the order 3 automorphisms. In particular, we consider two different aspects: the construction of isogenies of K3 surfaces which are not quotients map and the Shioda--Inose structures, which are relations between Abelian surfaces and K3 surfaces both admitting an order 3 symplectic automorphism.

\subsection{Isogenies} In \cite[Theorem 3.9]{CG} it is proved the following result
\begin{proposition}\label{prop: isometries lattice} There exists the following lattice isometry $$\left(\langle 6d\rangle\oplus K_{12}\right)'\stackrel{\simeq}{\longrightarrow} \langle 6d\rangle\oplus M_{\Z/3\Z}. $$
\end{proposition}
Let us denote by $X_d$ a K3 surface such that $\NS(X_d)\simeq \left(\langle 6d\rangle\oplus K_{12}\right)'$. 

By the previous proposition $\NS(X_d)\simeq \left(\langle 6d\rangle\oplus K_{12}\right)'\simeq \langle 6d\rangle\oplus M_{\Z/3\Z}$ and then:\\
$\bullet$ $X_d$ admits a symplectic automorphism $\sigma_d$ of order 3 (by Theorem \ref{theor:NS X}),\\
$\bullet$ there exists a K3 surface $S$ which admits a symplectic automorphism $\sigma_{S}$ of order 3 such that $X_{d}$ is the minimal resolution of $S/\sigma_{S}$ (by Theorem \ref{theor:NS Y}).

\begin{corollary}
There exists an infinite tower of isogenies
$$X_d\stackrel{\alpha_1}{\dashleftarrow}X_{3d}\stackrel{\alpha_2}{\dashleftarrow}X_{3^2d}\stackrel{\alpha_3}{\dashleftarrow}X_{3^3d}\ldots X_{3^{h-1}d}\stackrel{\alpha_h}{\dashleftarrow}X_{3^hd}\ldots.$$	
Each map $\alpha_i$ is the composition of a quotient map of order 3 and a birational map which resolves the singularities of the quotient. The compositions $\alpha_i\circ\alpha_{i-1}\circ \ldots\circ \alpha_{i-k}$ are not induced by quotients.
\end{corollary}
\proof To construct the tower, we start with a K3 surface $X_d$ such that $\NS(X_d)\simeq \langle 6d\rangle\oplus M_{\Z/3\Z}$. Since $M_{\Z/3\Z}$ is primitively embedded in $\NS(X_d)$, there exists a K3 surface $S$ such that $S$ admits a symplectic automorphism $\sigma$ and the minimal resolution of $S/\sigma$ is $X_d$. By Theorem \ref{theo: NS(Y) iff NS(X)}, the N\'eron--Severi group of $S$ is $\NS(S)\simeq \left(\langle 18d\rangle\oplus K_{12}\right)'\simeq \langle 18d\rangle\oplus M_{\Z/3\Z}$, hence $S$ is the surface that we denote by $X_{3d}$. Since the lattice $M_{\Z/3\Z}$ is primitively embedded in $\NS(X_{3d})$, there exists a K3 surface which admits a symplectic automorphism of order 3 such that the minimal resolution of the quotient of the surface by the automorphism is $X_{3d}$. This surface is $X_{3^2d}$, since its N\'eron--Severi group is isometric to $\left(\langle 54d\rangle\oplus K_{12}\right)'\simeq \langle 54\rangle\oplus M_{\Z/3\Z}$. This process can be iterated infinitely many times. 
As explained in \cite[Proposition 3.11]{CG}, it is not possible that the maps $X_h\dashrightarrow X_k$ are induced by quotient maps if $k/h>3$.\endproof

\subsection{Generalized Shioda--Inose structure}

The results in this paper are very much inspired by \cite{Mo}. In this paper Morrison also defined the Shioda--Inose structures, which are special relation between a K3 surface with a special symplectic involution and an Abelian surface. Morrison identified the K3 surfaces which admit a Shioda--Inose structure in \cite[Theorem 6.3]{Mo} by using the description he gave of the isometry induced by a symplectic involution on the second cohomology group of the K3 surfaces. Thanks to our description of the isometry induced by an order 3 symplectic automorphism on the second cohomology group of the K3 surfaces, one can generalize the results by Morrison. The details of this generalization are presented in the forthcoming paper \cite{GP}, but we recall here the main definition and theorem, which are the analogous the ones presented in \cite{Mo}, because they follow mainly by Section \ref{subsect: the action on Lambda} and in particular by Proposition \ref{prop: pi_*}.

We first generalize the definition of Shioda-Inose structure. The present definition is an extension of the definition given by \cite{OS}, see also \cite{CO}.
\begin{definition} A generalized Shioda--Inose structure of order 3 is given by $(A,\sigma_A, X,\sigma_X)$ where:
	\begin{itemize}\item $A$ is a 2-dimensional torus admitting a symplectic automorphism $\sigma_A$ of order 3 such that the minimal resolution of $A/\langle\sigma_A\rangle$ is a K3 surface $Km_3(A)$; \item $X$ is a K3 surface admitting a symplectic automorphism $\sigma_X$ of order 3 such that the minimal resolution of $X/\langle\sigma_X\rangle$ is the K3 surface $Km_3(A)$; \item $T_A\simeq T_X$.\end{itemize}\end{definition}
In particular a generalized Shioda--Inose structure of order 3 is associated to the following diagram
$$\xymatrix{A\ar@{-->}[dr]&&X\ar@{-->}[dl]\\&Km_3(A)&}$$
where the dash-arrows correspond to $3:1$ rational maps, and $T_A\simeq T_X$. 
	
The following theorem is the analogous of \cite[Theorem 6.3]{Mo} for the order 3 automorphisms and the proof is . Also the proof is very similar and detailed in \cite{GP}. It is based on the above results on order 3 symplectic automorphisms on K3 surfaces and on previous fundamental results on order 3 automorphisms on Abelian surfaces, due to \cite{F} and \cite{Ba2}.

	\begin{theorem}{\rm(\cite{GP})}\label{theo: generaliza Shioda--Inose} Let $X$ be a projective K3 surface. The following conditions are equivalent:
	\begin{itemize}
		\item[a)] There exists an Abelian surface $A$, an order 3 symplectic automorphism $\sigma_A$ on $A$, and an order 3 symplectic automorphism $\sigma_X$ on $X$ such that $(A,\sigma_A, X,\sigma_X)$ is a generalized Shioda--Inose structure of order 3.
		\item[b)] The transcendental lattice of $X$ is isometric to the one of an Abelian surface $A$ which admits an order 3 symplectic automorphism $\sigma_A$ such that $A/\langle\sigma_A\rangle$ is birationally equivalent to a K3 surface.
		\item[c)] There exists a primitive embedding $T_X\hookrightarrow U\oplus A_2(-1)$. 
		\item[d)] There exists a primitive embedding $(E_6\oplus E_6\oplus E_6)'\hookrightarrow NS(X)$.

		\end{itemize}
\end{theorem}

\end{document}